\documentclass[12pt]{amsart}
\usepackage{graphicx}
\usepackage{amssymb}
\usepackage{amsfonts}
\usepackage{amsmath}
\usepackage{array}
\usepackage{rotating}
\usepackage{bm}
\usepackage[matrix,frame]{xypic}

\newtheorem{example}{Example}[section]
\newtheorem{note}[example]{Note}
\newtheorem{theorem}[example]{Theorem}

\newtheorem{conjecture}[example]{Conjecture}
\newtheorem{definition}[example]{Definition}
\newtheorem{proposition}[example]{Proposition}

\newtheorem{lemma}[example]{Lemma}

\def\Proof{\noindent \it Proof -- \rm}
\def\qed{\hspace{3.5mm} \hfill \vbox{\hrule height 3pt depth 2 pt width 2mm}
\bigskip}

\def\id{{\rm id}}
\def\bp{BP}
\def\match{{\rm Ma}}

\def\sig{{\bm \sigma}}

\def\gf#1#2{\genfrac{}{}{0pt}{}{#1}{#2}}

\def\FQSym{{\bf FQSym}}
\def\MQSym{{\bf MQSym}}
\def\PQSym{{\bf PQSym}}
\def\CQSym{{\bf CQSym}}

\def\Sym{{\bf Sym}}
\def\NCSF{{\bf Sym}}
\def\QSym{{\it QSym}}
\def\FSym{{\bf FSym}}
\def\PBT{{\bf PBT}}
\def\sym{{\it Sym}}

\def\ssh{\Cup}

\def\sconc{\bullet}
\def\Std{{\rm Std}}
\def\cstd{{\bf Std}}
\def\Park{{\rm Park}}

\def\convol{{*}}

\def\park{{\bf a}}

\def\<{\langle}
\def\>{\rangle}

\def\C{\operatorname{\mathbb C}}
\def\K{\operatorname{\mathbb K}}
\def\Z{\operatorname{\mathbb Z}}
\def\Q{\operatorname{\mathbb Q}}
\def\N{\operatorname{\mathbb N}}

\def\F{{\bf F}}

\def\T{{\bf T}}
\def\G{{\bf G}}

\def\V{{\bf V}}

\def\SG{{\mathfrak S}}

\def\A{{\sf A}}

\def\Des{\operatorname{Des}}

\def\dim{{\rm dim}}

\def\PF{{\rm PF}}
\def\PPF{{\rm PPF}}

\def\shuff#1#2{\mathbin{
\hbox{\vbox{ \hbox{\vrule \hskip#2 \vrule height#1 width 0pt
}%
\hrule}%
\vbox{ \hbox{\vrule \hskip#2 \vrule height#1 width 0pt
\vrule }%
\hrule}%
}}}

\def\shuf{{\mathchoice{\shuff{7pt}{3.5pt}}%
{\shuff{6pt}{3pt}}%
{\shuff{4pt}{2pt}}%
{\shuff{3pt}{1.5pt}}}}%
\def\shuffle{\,\shuf\,}

\def\bij{{\phi}}
\newcommand{\free}[1]{\langle#1\rangle}
\def\conn{{\mathcal C}}
\def\connp{{\mathcal CP}}
\def\L{{\mathcal L}}
\def\LP{{\mathcal L'}}
\def\npn{{\bf n}}

\def\finer{\leq}
\def\MR{{\rm MR}}

\def\gf#1#2{\genfrac{}{}{0pt}{}{#1}{#2}}
\def\Ev{{\rm Ev}}
\def\intp{*}
\def\park{{\bf a}}

\def\A{{\bm A}}
\def\B{{\bm B}}
\def\X{{\mathbb X}}
\def\Y{{\mathbb Y}}
\def\T{{\mathbb T}}

\def\TD{{\mathfrak TD}}
\def\AA{{\mathcal A}} 
\def\KK{{\bf K}}

\def\qq{{\bf q}}
\def\maj{{\rm maj}}
\def\MAJ{{\rm MAJ}}

\def\Pp{{\bf P}}        
\def\pp{{\mathcal P}}   
\def\shape{{\rm shape\,}} 
\def\NCPQSym{{\bf NCPQSym}} 

\title[Free quasi-symmetric functions for wreath products]%
{Free quasi-symmetric functions and descent algebras for wreath products,
and noncommutative multi-symmetric functions}

\author[J.-C.~Novelli and J.-Y.~Thibon]%
{Jean-Christophe Novelli and Jean-Yves Thibon}

\address[] {Universit\'e Paris-Est, Institut Gaspard Monge  \\
5 Boulevard Descartes \\Champs-sur-Marne \\77454 Marne-la-Vall\'ee cedex 2 \\
FRANCE}
\email[Jean-Christophe Novelli]{novelli@univ-mlv.fr}
\email[Jean-Yves Thibon]{jyt@univ-mlv.fr} 
\date{\today}

\begin{document}

\begin{abstract}
We introduce analogs of the Hopf algebra of Free quasi-symmetric functions
with bases labelled by colored permutations. When the color set is a semigroup,
an internal product can be introduced.
This leads to the construction of  generalized descent algebras associated
with wreath products $\Gamma\wr\SG_n$ and to the corresponding generalizations
of quasi-symmetric functions. The associated Hopf algebras appear as natural
analogs of McMahon's multisymmetric functions. As a consequence, we obtain an
internal product on ordinary multi-symmetric functions.
We extend these constructions to Hopf algebras of colored parking
functions, colored non-crossing partitions and parking functions of type $B$.
\end{abstract}

\maketitle
\tableofcontents

\section{Introduction}

The Hopf algebra of Free Quasi-Symmetric Functions $\FQSym$ \cite{NCSF6}
is an algebra of noncommutative polynomials associated with the
sequence $(\SG_n)_{n\geq0}$ of all symmetric groups.
It is connected by Hopf homomorphisms to several other important algebras
associated with the same sequence of groups : Free symmetric functions (or
coplactic algebra) $\FSym$ \cite{PR,NCSF6}, Non-commutative symmetric
functions (or descent algebras) $\NCSF$ \cite{NCSF1}, Quasi-Symmetric
functions $\QSym$ \cite{Ge84}, Symmetric functions $\sym$, and also, Planar
binary trees $\PBT$ \cite{LR1,HNT}, Matrix quasi-symmetric functions $\MQSym$
\cite{NCSF6,Hiv}, Parking functions $\PQSym$ \cite{KW,NT}, and so on.

Most of these Hopf algebras are endowed with an internal product, generalizing
the one of ordinary symmetric functions. The basic example is provided by
noncommutative symmetric functions, whose homogeneous components
can be identified with the Solomon descent algebras of symmetric groups
\cite{NCSF1}.

Symmetric groups are the Coxeter groups of type $A$, and there are
descent algebras for other types as well. However, the direct sums of
the descent algebras of types $B$ or $D$ are not Hopf algebras in any natural
way. But there are Hopf algebras associated with wreath products
$\Z_\ell \wr \SG_n$, the Mantaci-Reutenauer algebras \cite{MR}, which admit
internal products, and contain the Solomon algebras of type $B$ for $\ell=2$.

From the point of view of symmetric functions, $\MR^{(\ell)}$, the Mantaci-Reutenauer algebra
of level $\ell$ is the free product of $\ell$ copies of $\Sym$.
It is therefore the natural noncommutative analog of 
$(\sym)^{\otimes \ell} \simeq \sym(X_0 ; \ldots ; X_{\ell-1})$, the algebra
of symmetric functions in $\ell$ independent alphabets, which is also the
Grothendieck ring of the tower of algebras $(\C[\Z_\ell \wr \SG_n])$.
And indeed, it has been shown in 
\cite{HNT-AKS} that $\MR^{(\ell)}$ was the Grothedieck ring of projective modules
over the 0-Ariki-Koike-Shoji algebras, a degeneracy of the Hecke algebras
associated with $\Z_\ell \wr \SG_n$.

However, with $\ell$ independent alphabets, one can build a larger Hopf
algebras. In the commutative case, it is the algebra of {\em multi-symmetric
functions}, first introduced by McMahon \cite{McM}, and briefly investigated
by Gessel from a modern point of view in \cite{Ge87}. It is defined as follows.
Setting $X_i=\{x_{i,j}|j=1,\ldots,n\}$, the multi-symmetric polynomials
are the invariants of $\SG_n$ in $\C[X_0,\ldots,X_{\ell-1}]$ for the diagonal
action (by the automorphisms $\sigma(x_{i,j})=x_{i,\sigma(j)}$).
This is an algebra, which, as usual, acquires a Hopf algebra structure
in the limit $n\rightarrow\infty$.

In the following, we will start with a level $\ell$ analogue of $\FQSym$, 
whose bases are
labelled by $\ell$-colored permutations. Imitating the embedding of $\NCSF$ in
$\FQSym$, we obtain a Hopf subalgebra of level $\ell$ called $\NCSF^{(\ell)}$,
which is a natural noncommutative analog of McMahon's algebra of multi-symmetric functions,
and turns out to be dual to Poirier's quasi-symmetric functions. Its
homogenous components can be endowed with an internal product, thus providing an
analog of Solomon's descent algebras for wreath products, bigger than the 
Mantaci-Reutenauer algebras, and in which most useful properties such as
the splitting formula remain valid. By commutative image, this yields
an internal product on multi-symmetric functions.

The Mantaci-Reutenauer descent algebra $\MR^{(\ell)}$ arises as a natural Hopf
subalgebra of $\NCSF^{(\ell)}$ and its dual is computed in a straightforward
way by means of an appropriate Cauchy formula.

Finally, we introduce a Hopf algebra of colored parking functions,
and use it to define Hopf algebras structures on parking functions
and non-crossing partitions of type $B$.

The main results of this paper have been announced in the draft
\cite{NT-level}. Since then, some of these results, in particular the construction
of $\NCSF^{(\ell)}$, have been used by Baumann and Hohlweg \cite{BH}, whose
paper provide detailed proofs. Hence, we shall only include the proofs which 
cannot be found in their paper. In particular, we propose an alternative
approach to the internal product, which is introduced by a duality argument,
and derive its main properties from those of the dual coproduct.

{\footnotesize
{\it Acknowledgements.}
This project has been partially supported by the grant ANR-06-BLAN-0380.
The authors would also like to thank the contributors of the MuPAD project,
and especially of the combinat part, for providing the development environment
for their research (see~\cite{HTm} for an introduction to MuPAD-Combinat).
}
\section{Background and notations}

We first explain how to adapt the classical definitions and operations to the
$\ell$-colored case.

\subsection{Colored alphabets}

We shall start with an $\ell$-colored alphabet
\begin{equation}
\A = A^0 \sqcup A^1 \sqcup \cdots \sqcup A^{\ell-1},
\end{equation}
such that all $A^i$ are of the same cardinality $N$, which will be assumed to
be infinite in the sequel.
Let $C$ be the alphabet $\{c_0,\ldots,c_{\ell-1}\}$ and $A$ be the auxiliary
ordered alphabet $\{1,2,\ldots\}$ (the letter $C$ stands for \emph{colors}
and $A$ for \emph{alphabet}) so that $\A$ can be identified with the cartesian
product $A\times C$:
\begin{equation}
\A \simeq A \times C = \{ (a,c), a\in A,\ c\in C \}.
\end{equation}
A colored letter $(i,c)$ will be denoted in bold type ${\bm i}$.
Given two colored words, their \emph{concatenation} is obtained by
concatenating separately the elements coming from $A$ and from $C$.
We will sometimes allow $\ell=\infty$.

\subsection{Colored standardization}

Let $w$ be a word in $\A$, represented as $(v,u)$ 
with $v\in A^*$ and $u\in C^*$.
Then the \emph{colored standardized word} $\cstd(w)$ of $w$ is
\begin{equation}
\cstd(w) := (\Std(v),u), 
\end{equation}
where $\Std(v)$ is the usual standardization on words.

Recall that the standardization process sends a word $v$ of length $n$ to a
permutation $\Std(v)\in\SG_n$, called the \emph{standardized} of $v$, defined
as the permutation obtained by iteratively scanning $v$ from left to right,
and labelling $1,2,\ldots$ the occurrences of its smallest letter, then
numbering the occurrences of the next one, and so on. Alternatively, $\Std(v)$
is the permutation having the same inversions as $v$.

For example, $\Std({abcadbdaa})=157286934$:
\begin{equation}
\begin{array}{ccccccccc}
a  & b  & c  & a  & d  & b  & d  & a  & a\\
a_1& b_5& c_7& a_2& d_8& b_6& d_9& a_3& a_4\\
1   &5   &7   &2   &8   &6   &9   &3   &4
\end{array}
\end{equation}
so that
\begin{equation}
\cstd(abcadbdaa,144120100) = (157286934,144120100).
\end{equation}

\subsection{Colored shifted operations}

For an element $v=(v_1,v_2,\ldots)$ of $A$ and an integer $k$, denote by
$v[k]$ the \emph{shifted word} $(v_1+k)\cdots (v_n+k)$, \emph{e.g.},
$312[4]=756$.
Given a colored word
${\bm \alpha} = (\alpha,u)$, we set
${\bm \alpha}[k]=(\alpha[k],u)$.

The \emph{shifted concatenation} of two words $v$ and $v'$ is defined by
\begin{equation}
v\bullet v' := v\cdot v'[k]
\end{equation}
where $k$ is the length of $v$.

The \emph{shifted concatenation of two colored words} $(v,c)$, $(v',c')$
is defined by
\begin{equation}
(v,c)\bullet (v',c') := (v\cdot v'[k], c\cdot c').
\end{equation}
For example,
\begin{equation}
(13241,00322)\bullet(12,23) = (1324167,0032223).
\end{equation}

Finally, recall that the \emph{shuffle product} of two words $au$ and $bv$ is
defined by
\begin{equation}
au\shuffle bv = a (u\shuffle bv) + b (au\shuffle v),
\end{equation}
where $a$ and $b$ are letters and $u$ and $v$ are words,
with the initial conditions
\begin{equation}
u\shuffle \epsilon = \epsilon \shuffle u = u,
\text{\qquad\qquad $\epsilon$ being the empty word}.
\end{equation}

\noindent
This extends to the colored case, considering colored
words as the concatenation of biletters.

The \emph{shifted shuffle product} is
\begin{equation}
u\ssh v := u\shuffle v[k],
\end{equation}
where $k$ is the size of $u$.

\section{Free quasi-symmetric functions of level $\ell$}

\subsection{$\FQSym^{(\ell)}$ and $\FQSym^{(\Gamma)}$}

A \emph{colored permutation} is a pair $(\sigma,u)$, with $\sigma\in\SG_n$ and
$u\in C^n$, the integer $n$ being the \emph{size} of this permutation.

\begin{definition}
The \emph{dual free colored quasi-ribbon} $\G_{\sigma,u}$ labelled by a
colored permutation $(\sigma,u)$ of size $n$ is the noncommutative
polynomial
\begin{equation}
\label{defG}
\G_{\sigma,u} := \sum_{w\in A^n ; \cstd(w)=(\sigma,u)} w \quad\in\Z\free{\A}.
\end{equation}
\end{definition}

Recall that the \emph{convolution} of two permutations $\sigma$ and $\mu$ is
the set $\sigma\convol\mu$ (identified with the formal sum of its elements)
of permutations $\tau$ such that the standardized word of the $|\sigma|$ first
letters of $\tau$ is $\sigma$ and the standardized word of the remaining
letters of $\tau$ is $\mu$ (see~\cite{Reu}).
We then have:

\begin{theorem}
\label{prodG}
Let $(\sigma',u')$ and $(\sigma'',u'')$ be colored permutations.
Then
\begin{equation}
\G_{\sigma',u'}\,\,\G_{\sigma'',u''} = \sum_{\sigma\in \sigma'\convol\sigma''}
\G_{\sigma,u'\cdot u''}.
\end{equation}
Therefore, the dual free colored quasi-ribbons span a $\Z$-subalgebra 
${\FQSym^{(\ell)}}$ of the free associative algebra.
\end{theorem}

\Proof
This is immediate from the product of the usual free quasi-symmetric
functions $\G_\sigma$:
\begin{equation}
\label{prodG-eq}
\G_{\sigma'}\,\,\G_{\sigma''}
= \sum_{\sigma\in \sigma'\convol\sigma''}  \G_{\sigma}.
\end{equation}
\qed

Note that all colored permutations indexing a product of $\G$ have given
colors at the same places.
For example,
\begin{equation}
\label{prodGex}
\begin{split}
\G_{(21,41)}\G_{(12,31)}
= &
\ \ \ \G_{(2134,4131)} + \G_{(3124,4131)} + \G_{(4123,4131)} \\
  &+  \G_{(3214,4131)} + \G_{(4213,4131)} + \G_{(4312,4131)}.
\end{split}
\end{equation}

One can define a coproduct  by the usual trick of
sums of alphabets: observe that we only need a total order on $A$ to define
the colored standardisation, so that taking two isomorphic copies $A'$ and
$A''$ of $A$, we define $\A'\oplus \A''$ as $(A'\oplus A'')\times C$, where
$A'\oplus A''$ denotes the ordered sum. Assuming furthermore that $\A'$ and
$\A''$ commute, we identify $f(\A') g(\A'')$ with $f\otimes g$ and define a
coproduct by:
\begin{equation}
\Delta f(\A) = f(\A'\oplus \A'').
\end{equation}

By construction, this is an algebra morphism from
${\FQSym^{(\ell)}}$ to ${\FQSym^{(\ell)}} \otimes {\FQSym^{(\ell)}}$,
so that

\begin{theorem}
${\FQSym^{(\ell)}}$ is a graded connected bialgebra.
Hence, it is a Hopf algebra.
The coproduct is given by
\begin{equation}
\label{deltaG}
\Delta \G_{\sigma,u} :=
 \sum_{\gf{(\sigma',\sigma'',u',u'')}%
          {(\sigma,u)\in(\sigma',u')\ssh(\sigma'',u'')}}
 \G_{(\sigma',u')} \otimes \G_{(\sigma'',u'')}.
\end{equation}
\end{theorem}

For example,
\begin{equation}
\label{exDelG}
\begin{split}
\Delta \G_{3142,2412} =
&\ 1\otimes \G_{3142,2412} + \G_{1,4}\otimes \G_{231,212} +
   \G_{12,42}\otimes \G_{12,21} \\
& + \G_{312,242}\otimes \G_{1,1} + \G_{3142,2412}\otimes 1.
\end{split}
\end{equation}

\Proof
This is again immediate from the coproduct of the usual free
quasi-symmetric functions $\G_\sigma$:
\begin{equation}
\Delta \G_{\sigma}
= \sum_{\sigma\in \sigma'\ssh\,\sigma''}  \G_{\sigma'}\otimes \G_{\sigma''}.
\end{equation}
\qed

\subsection{Duality in $\FQSym^{(\ell)}$}

\begin{definition}
The \emph{free $\ell$-quasi-ribbon} $\F_{\sigma,u}$ labelled by a colored
permutation $(\sigma,u)$ is the noncommutative polynomial
\begin{equation}
\label{defF}
\F_{\sigma,u} := \G_{\sigma^{-1},u\cdot\sigma^{-1}},
\end{equation}
where the action of a permutation on the right of a word permutes the
positions of the letters of the word.
\end{definition}

For example,
\begin{equation}
\F_{3142,2142} = \G_{2413,1422}\,.
\end{equation}

The product and coproduct of the $\F_{\sigma,u}$ can be easily described in
terms of shifted shuffle and deconcatenation of colored permutations.

\begin{theorem}
\label{prodcopF}
Let $\sig'$ and $\sig''$ be two colored permutations.
Then
\begin{equation}
\label{prodF}
\F_{\bm \sigma'} \F_{\bm \sigma''} =
\sum_{\sig \in \sig'\ssh\sig''} \F_{\sig},
\end{equation}

\noindent
and
\begin{equation}
\label{copF}
\Delta\F_{\sig} =
\sum_{\gf{w',w''}{\sig=w'.w''}}
 \F_{\cstd(w')}\otimes \F_{\cstd(w'')}.
\end{equation}
\end{theorem}

\Proof
Without colors, these formulas are the usual product of coproduct formulas of
the $\F$ in $\FQSym$. With colors, one just has to observe that colors follow
the letter to which they are attached.
\qed

Note that all colored permutations indexing a product of $\F$ have given
colors associated with the same values, which is consistent with the
corresponding remark on the $\G$ since places and values are exchanged when
taking the inverse of a permutation.

For example, compare the following with Equation~(\ref{prodGex}):
\begin{equation}
\begin{split}
\F_{(21,14)}\F_{(12,31)}
= &
\ \ \ \F_{(2134,1431)} + \F_{(2314,1341)} + \F_{(2341,1314)} \\
  &+  \F_{(3214,3141)} + \F_{(3241,3114)} + \F_{(3421,3114)}.
\end{split}
\end{equation}

Here is an example of coproduct on the $\F$ basis:
\begin{equation}
\begin{split}
\Delta \F_{(23514,14212)}
=&\ 1 \otimes \F_{(23514,14212)} +
   \F_{(1,1)} \otimes \F_{(2413,4212)} +
   \F_{(12,14)} \otimes \F_{(312,212)}\\
&\!\!+ \F_{(123,142)} \otimes \F_{(12,12)} +
   \F_{(2341,1421)} \otimes \F_{(1,12)} +
   \F_{(23514,14212)} \otimes 1.
\end{split}
\end{equation}

Let us define a scalar product on $\FQSym^{(\ell)}$ by
\begin{equation}
\langle \F_{\sigma,u} , \G_{\sigma',u'} \rangle :=
  \delta_{ \sigma,\sigma'} \delta_{u,u'},
\end{equation}
where $\delta$ is the Kronecker symbol.

\begin{theorem}
For any $U,V,W\in\FQSym^{(\ell)}$,
\begin{equation}
\langle \Delta U, V\otimes W \rangle =
\langle U, V W \rangle,
\end{equation}
so that, $\FQSym^{(\ell)}$ is self-dual: the map
$\F_{\sigma,u} \mapsto {\G_{\sigma,u}}^*$ is an isomorphism from
$\FQSym^{(\ell)}$ to its graded dual.
\end{theorem}

\Proof
Straightforward from Theorem~\ref{prodcopF}.
\qed

\begin{note}
{\rm
Let $\bij$ be any bijection from $C$ to $C$, extended to words by
concatenation. Then if one defines the free $\ell$-quasi-ribbon as
\begin{equation}
\F_{\sigma,u} := \G_{\sigma^{-1},\bij(u)\cdot\sigma^{-1}},
\end{equation}
the previous theorems remain valid since one only permutes the labels of
the basis $(\F_{\sigma,u})$.
Moreover, if $C$ has a semigroup structure, the colored permutations
$(\sigma,u)\in\SG_n\times C^n$ can be interpreted as elements of the
semi-direct product $H_n := \SG_n\ltimes C^n$ with multiplication rule
\begin{equation}
\label{wreath}
(\sigma ; c_1,\ldots,c_n) \cdot (\tau ; d_1,\ldots,d_n) :=
(\sigma\tau ; c_{\tau(1)}d_1, \ldots, c_{\tau(n)}d_n).
\end{equation}
In furthermore $C$ is a group, one can choose $\bij(\gamma):=\gamma^{-1}$ and define the scalar
product as before, so that the adjoint basis of the $(\G_{h})$ becomes
$\F_h := \G_{h^{-1}}$.
In the sequel, we will be mainly interested in the cases $C:=\Z/\ell\Z$, and we
will indeed make that choice for $\bij$ whenever $C$ is a group.
}
\end{note}

\subsection{Algebraic structure of $\FQSym^{(\ell)}$}

Recall that a permutation $\sigma$ of size $n$ is \emph{connected}
\cite{MR,NCSF6} if, for any $i<n$, the set $\{\sigma_1,\ldots,\sigma_i\}$ is
different from $\{1,\ldots,i\}$.

We denote by $\conn$ the set of connected permutations, and by
$c_n:=|\conn_n|$ the number of such permutations in $\SG_n$. For later
reference, we recall that the generating series of $c_n$ is Sequence~A003319
of~\cite{Slo}:
\begin{equation}
c(t) := \sum_{n\ge 1} c_n t^n
  = 1 - \left(\sum_{n\ge 0} n! t^n\right)^{-1}\\
  = t+{t}^{2}+3\,{t}^{3}+13\,{t}^{4}+71\,{t}^{5}+461\,{t}^{6} +O(t^7)\,.
\end{equation}

Let the \emph{connected colored permutations} be the $(\sigma,u)$ with
$\sigma$ connected and $u$ arbitrary. Their generating series is given by
$c(\ell t)$.

From~\cite{NCSF6}, we immediately get
\begin{proposition}
$\FQSym^{(\ell)}$ is free over the set $\F_{\sigma,u}$ (or $\G_{\sigma,u}$),
where $(\sigma,u)$ is connected.\qed
\end{proposition}

For example, the generating series of the algebraic generators of
$\FQSym^{(2)}$ is
\begin{equation}
2\,t + 4\,t^2 + 24\,t^3 + 208\,t^4 + 2272\,t^5 + 29504\,t^6 + 441216\,t^7
+\dots
\end{equation}

\subsection{Primitive elements of $\FQSym^{(\ell)}$}

Let $\L^{(\ell)}$ be the primitive Lie algebra of $\FQSym^{(\ell)}$.
Since $\Delta$ is not cocommutative, $\FQSym^{(\ell)}$ cannot be the universal
enveloping algebra of $\L^{(\ell)}$.
But since it is cofree, it is, according to~\cite{LRdip}, the universal
enveloping dipterous algebra of its primitive part $\L^{(\ell)}$.

Let $\G^{\sigma,u}$ be the multiplicative basis defined by
$\G^{\sigma,u}=\G_{\sigma_1,u_1}\cdots\G_{\sigma_r,u_r}$ where
$(\sigma,u)=(\sigma_1,u_1)\bullet\cdots\bullet(\sigma_r,u_r)$ is the unique
maximal factorization of $(\sigma,u)\in\SG_n\times C^n$ into connected colored
permutations.

\begin{proposition}
Let $\V_{\sigma,u}$ be the adjoint basis of $\G^{\sigma,u}$.
Then, the family $(\V_{\alpha,u})_{\alpha\in\conn}$ is a basis of
$\L^{(\ell)}$.
In particular, we have $\dim\,\L^{(\ell)}_n =\ell^n c_n$.
%
Moreover, $\L^{(\ell)}$ is free.
\end{proposition}

\Proof
The first part of the statement follows from~\cite{NCSF6}.
The second part comes from the fact that $\FQSym^{(\ell)}$ is bidendriform
(Theorem~\ref{fqsl-bidend} below).
\qed

For example, since $\L^{(\ell)}$ is free, the generating series by degree
of its generators is (with $\ell=2$):
\begin{equation}
  \begin{split}
    1 - \prod_{n\geq1}{(1-t^n)}^{\ell^n c_n}
    &= 1-(1-t)^2(1-t^2)^4(1-t^3)^{24} \cdots\\
    &= 2\,t + 3\,t^2 + 16\,t^3 + 158\,t^4 + 1\,796\,t^5\\
    & + 24\,250\,t^6 + 372\,656\,t^7 + \dots\\
  \end{split}
\end{equation}

\noindent
and the Hilbert series of the universal enveloping algebra of
$\FQSym^{(\ell)}$ (its domain of cocommutativity) is, again with $\ell=2$,
\begin{equation}
  \begin{split}
    \prod_{n\geq1}{(1-t^n)}^{-\ell^n c_n}
    &= 1 + 2\,t + 7\,t^2 + 36\,t^3 + 283\,t^4 + 2\,898\,t^5 + 36\,169\,t^6\\
      &+ 524\,976\,t^7 + \dots \\
  \end{split}
\end{equation}

\subsection{Dendriform structure of $\FQSym^{(\ell)}$}

Foissy introduced the notion of \emph{bidendriform} bialgebras~\cite{Foi},
generalizing the notion of \emph{dendrifom}
algebras~(\emph{cf.} \cite{Lod}) and proved some conjectures about $\FQSym$,
presented in~\cite{NCSF6}. We shall adapt this technology to the colored case.
We shall not recall all the theory, since complete details can be found
in~\cite{Foi}.

Recall that the generators of $\FQSym^{(\ell)}$ as a dendrifom algebra are
called \emph{totally primitive elements}~\cite{Foi} and that their generating
series is given by
\begin{equation}
{\rm TP} := \frac{{\rm PQ}-1}{{\rm PQ}^2},
\end{equation}
where ${\rm PQ}$ is the generating series of $\FQSym^{(\ell)}$.

\begin{theorem}
\label{fqsl-bidend}
The algebra $\FQSym^{(\ell)}$ has a structure of bidendriform
bialgebra~\cite{Foi}, hence is free as a Hopf algebra and as a dendriform
algebra, cofree, self-dual, and its primitive Lie algebra is free.

Moreover, the totally primitive elements of $\FQSym^{(\ell)}$ are the totally
primitive elements of $\FQSym$ with any coloring.
\end{theorem}

\Proof
It has been done by Foissy in~\cite{Foi} in the case of $\FQSym$.
But since the dendriform and codendriform structure do not involve the color
alphabet $C$, the property is true in this case as well.

Similarly, colors do not play any role in determining if a given element
is (totally) primitive or not.
\qed

For example, the dendriform generators of $\FQSym$ have as degree generating
series
\begin{equation}
\sum_i {\rm dg}_i t^i = t + t^3 + 6\,t^4 + 39\,t^5 + 284\,t^6 + 2305\,t^7
+ \dots
\end{equation}
so that the dendriform generators of $\FQSym^{(2)}$ have as degree generating
series $2^i{\rm dg}_i$:
\begin{equation}
2\,t + 8\,t^3 + 96\,t^4 + 1248\,t^5 + 18176\,t^6 + 295040\,t^7
+\dots,
\end{equation}

Note that there cannot be any relation, even dendriform relations, among the
elements $\F_{1,c}$ where $c\in C$, so that $\FQSym^{(l)}$ contains the free
dendriform algebra $\PBT^{(\ell)}$ on $\ell$ generators.

\subsection{Internal product of $\FQSym^{(\ell)}$}

When $C$ is a semigroup, an internal product can be defined on
$\FQSym^{(\ell)}$ by
\begin{equation}
\F_{\sigma,u}*\F_{\tau,v}=\F_{\mu,w}
\end{equation}
where $(\mu,w)$ is the product $(\sigma,u)\cdot(\tau,v)$ in the
wreath product, defined by formula (\ref{wreath}), that is

\begin{equation}
\label{intpF}
\F_{(\sigma,u)} \intp \F_{(\tau,v)} = \F_{(\sigma\tau, (u\tau)\cdot v)},
\end{equation}
where $u\tau$ is the word $u_{\tau_{1}}\dots u_{\tau_{n}}$ and
$u\cdot v$ is the componentwise product defined by
$(u_1v_1,\dots,u_nv_n)$.

For example, if the color group is $\Z$
\begin{equation}
\F_{(1324,1011)} \intp \F_{(2413,3200)} =
\F_{(3412, 3311)}.
\end{equation}
\begin{equation}
\F_{(165324,102011)} \intp \F_{(625413,322011)} =
\F_{(462315, 423023)}.
\end{equation}
%
This can be reduced to any $\Z/\ell\Z$, \emph{e.g.}, with $\ell=3$,
\begin{equation}
\F_{(165324,102011)} \intp \F_{(625413,022011)} =
\F_{(462315, 120020)}.
\end{equation}

In the $\G$ basis, one has
\begin{equation}
\label{intpG}
\G_{(\sigma,u)} \intp \G_{(\tau,v)} = \G_{(\tau\sigma, u.(v\sigma))},
\end{equation}

\section{Noncommutative symmetric functions of level $\ell$}

\subsection{$\ell$-partite numbers}

Following McMahon~\cite{McM}, we define an \emph{$\ell$-partite number}
$\npn=(n_1,\dots,n_\ell)$
as a column vector in $\N^\ell$, and a \emph{vector composition of $\npn$} of
weight $|\npn|:=\sum_{1}^\ell{n_i}$ and length $m$ as a $\ell\times m$ matrix
$\bf I$ of nonnegative integers, with row sums vector $\npn$ and no zero
column.

For example,
\begin{equation}
\label{exM}
{\bf I} =
\begin{pmatrix}
1 & 0 & 2 & 1 \\
0 & 3 & 1 & 1 \\
4 & 2 & 1 & 3 \\
\end{pmatrix}
\end{equation}
is a vector composition (or a $3$-composition, for short)
of the $3$-partite number
\scriptsize$\begin{pmatrix} 4\\ 5\\ 10\end{pmatrix}$ \normalsize
of weight $19$ and length $4$.

For each $\npn\in\N^\ell$  of weight $|\npn|=n$, we define a
\emph{level $\ell$ complete homogeneous noncommutative symmetric function} as
\begin{equation}
\label{defS}
S_{\npn} := \sum_{u ; |u|_i=n_i} \G_{1\cdots n, u}.
\end{equation}
It is the sum of all possible colorings of the identity permutation with $n_i$
occurrences of color $i$ for each $i$.

\subsection{The Hopf algebra $\NCSF^{(\ell)}$}

Let $\NCSF^{(\ell)}$ be the subalgebra of $\FQSym^{(\ell)}$ generated by the
$S_{\npn}$ (with the convention $S_{\bf 0}=1$).
The Hilbert series of $\NCSF^{(\ell)}$ is easily found to be
\begin{equation}
S_\ell(t) :=
  \sum_{n} {\dim\ \NCSF_n^{(\ell)}t^n} = \frac{(1-t)^\ell}{2(1-t)^\ell-1}.
\end{equation}

For example, with $\ell=2$, one has
\begin{equation}
S_2(t) :=  1 + 2\,t + 7\,t^2 + 24\,t^3 + 82\,t^4 + 280\,t^5 + 956\,t^6
+ 3264\,t^7 + \dots
\end{equation}
which is Sequence~A003480 of~\cite{Slo}.

For general $\ell$, it is well-known in the combinatorial folklore (and easy to
prove by means of generating series expansions) that, for all
$n\geq1$,
\begin{equation}
{\rm ncs}_n(\ell) := n!\ \dim( \NCSF_n^{(\ell)})
             = \sum_{k=1}^n S(n,k)\, p_k\, \ell^k
\end{equation}
where $S(n,k)$ is the sequence of absolute values of Stirling numbers of the
first kind (sequence~A130534 of~\cite{Slo}) and $p_k$ is the sequence of
ordered Bell numbers (also known as packed words or preferential arrangements,
Sequence~A000670 of~\cite{Slo}).

\begin{theorem}
$\NCSF^{(\ell)}$  is free over the set $\{S_{\npn}, |\npn|>0 \}$, so that a
linear basis is given by
\begin{equation}
S^{\bf I} = S_{{\bf i}_1} \cdots S_{{\bf i}_m},
\end{equation}
where ${\bf i}_1,\cdots,{\bf i}_m$ are the columns of $\bf I$.

Moreover, $\NCSF^{(\ell)}$ is a Hopf subalgebra of $\FQSym^{(\ell)}$
and the coproduct of the generators is given by
\begin{equation}
\Delta S_\npn = \sum_{{\bf i}+{\bf j}= {\bf n}} S_{\bf i}\otimes S_{\bf j},
\end{equation}
where the sum ${\bf i}+{\bf j}$ is taken in the space $\N^l$. In particular,
$\NCSF^{(\ell)}$ is cocommutative.
\end{theorem}

\Proof
Consider a linear relation between the $S^{\bf I}$.
Since without colors the algebra is free (it is $\NCSF$), the linear relation
splits into many linear relations involving terms of the form
$S^{{\bf i}_1\dots{\bf i}_r}$ with $|{\bf i}_1|, \dots, |{\bf i}_r|$ fixed.
But there are no relations of this form, thanks to the definition of the
$\G_{\sig}$ as a sum of colored words.

Given the coproduct of the $\G_{\sig}$ of $\FQSym^{(\ell)}$ and since the
$S_{\bf n}$ are sums of $\G$, the coproduct of an $S_{\bf n}$ amounts to
unshuffling the color words, whence their coproduct formula.
\qed

For example,
\begin{equation}
\begin{split}
\Delta S^{\hbox{\scriptsize$\begin{pmatrix}1\\0\\2\end{pmatrix}$}}
=&\ \ \ \
S^{\hbox{\scriptsize$\begin{pmatrix}1\\0\\2\end{pmatrix}$}}\otimes
S^{\hbox{\scriptsize$\begin{pmatrix}0\\0\\0\end{pmatrix}$}}
+
S^{\hbox{\scriptsize$\begin{pmatrix}0\\0\\2\end{pmatrix}$}}\otimes
S^{\hbox{\scriptsize$\begin{pmatrix}1\\0\\0\end{pmatrix}$}}
+
S^{\hbox{\scriptsize$\begin{pmatrix}1\\0\\1\end{pmatrix}$}}\otimes
S^{\hbox{\scriptsize$\begin{pmatrix}0\\0\\1\end{pmatrix}$}} \\[7pt]
&+
S^{\hbox{\scriptsize$\begin{pmatrix}0\\0\\1\end{pmatrix}$}}\otimes
S^{\hbox{\scriptsize$\begin{pmatrix}1\\0\\1\end{pmatrix}$}} 
+
S^{\hbox{\scriptsize$\begin{pmatrix}1\\0\\0\end{pmatrix}$}}\otimes
S^{\hbox{\scriptsize$\begin{pmatrix}0\\0\\2\end{pmatrix}$}}
+
S^{\hbox{\scriptsize$\begin{pmatrix}0\\0\\0\end{pmatrix}$}}\otimes
S^{\hbox{\scriptsize$\begin{pmatrix}1\\0\\2\end{pmatrix}$}}
.
\end{split}
\end{equation}

\subsection{Algebraic structure of $\NCSF^{(\ell)}$}

The number of generators of $\NCSF^{(\ell)}$ of degree $n$ is
given by the number of $\ell$-partite numbers of total sum $n$.
So its generating series is
\begin{equation}
(1-t)^{-\ell}-1=\sum_{n\geq1} \binom{\ell+n-1}{n} t^n.
\end{equation}

We shall denote by $G(\ell)$ the set of 
nonzero $\ell$-partite numbers.

\subsection{Primitive elements of $\NCSF^{(\ell)}$}
\label{prim-NC}

$\NCSF^{(\ell)}$ being a graded connected cocommutative Hopf algebra, it
follows from the Cartier-Milnor-Moore theorem that it is the universal
enveloping algebra of $L^{(\ell)}$:
\begin{equation}
\NCSF^{(\ell)} = U(L^{(\ell)}),
\end{equation}
where $L^{(\ell)}$ is the Lie algebra of its primitive elements.
Let us now prove
\begin{theorem}
\label{prim-libNC}
As a graded Lie algebra, the primitive Lie algebra $L^{(\ell)}$
of $\NCSF^{(\ell)}$ is free over a set indexed by $G(\ell)$.
\end{theorem}

\Proof
If $L^{(\ell)}$ is free, by standard arguments on generating series, the
number of generators of $L^{(\ell)}$ in degree $n$ must be the number of
algebraic generators of $\NCSF^{(\ell)}$ in degree $n$, parametrized for
example by $G(\ell)$.
We will now show that $L^{(\ell)}$ has at least this number of generators and
that those generators are algebraically independent, determining completely
thee dimensions of the homogeneous components $L^{(\ell)}_n$ of
$L^{(\ell)}$ whose generating series begins by
\begin{equation}
t + \binom{\ell+1}{2} t^2
  + \left(\binom{\ell+1}{2}+\binom{\ell+2}{3} \right) t^3 + \dots
\end{equation}
Following Reutenauer~\cite{Reu} p.~58, denote by $\pi_1$ the Eulerian
idempotent, that is, the endomorphism of $\NCSF^{(\ell)}$ defined by
$\pi_1=\log^*(Id)$. It is obvious, thanks to the definition of $S_{\bf p}$
that
\begin{equation}
\pi_1(S_{\bf p}) = S_{\bf p} + \cdots,
\end{equation}
where the dots stand for terms $S_{\bf I}$ such that ${\bf I}$ are vector
compositions with strictly more than one column.
Since the $S_{\bf p}$ where ${\bf p}$ is a $\ell$-partite number are
algebraically independent, the dimension of $L^{(\ell)}_n$ is at least equal
to the cardinality of the elements of $G(l)$ of size $n$. So $L^{(\ell)}$ is
indeed free over a set of primitive elements parametrized by $G(\ell)$.
\qed

For example, with $\ell=2$, the generating series of the dimensions of
$L^{(\ell)}$ is
\begin{equation}
1+ 2\,t + 4\,t^2 + 12\,t^3 + 31\,t^4 + 92\,t^5 + 256\,t^6 +
772\,t^7 + \dots
\end{equation}
With $\ell=3$, one finds
\begin{equation}
1+ 3\,t + 9\,t^2 + 36\,t^3 + 132\,t^4 + 534\,t^5 + 2140\,t^6 +
8982\,t^7 + \dots
\end{equation}

\subsection{Duality in $\NCSF^{(\ell)}$}

Recall that the underlying colored alphabet $\A$ can be seen as
$A^0 \sqcup \cdots \sqcup A^{\ell-1}$, with $A^i = \{ a^{(i)}_j, j\geq1 \}$.
Let ${\bf x} = (x^{(0)}, \ldots, x^{(\ell-1)})$, where the $x^{(i)}$ are
$\ell$ commuting variables.
In terms of $\A$, the generating function of the complete functions can be
written as
\begin{equation}
\sigma_{\bf x}(\A) = \prod_{i\geq0}^{\rightarrow}
\left(1-\sum_{0\leq j\leq \ell-1} x^{(j)} a_{i}^{(j)} \right)^{-1}
= \sum_\npn {S_{\bf n}(\A) {\bf x}^{\bf n}},
\end{equation}
where ${\bf x}^{\bf n} = (x^{(0)})^{n_0} \cdots (x^{(\ell-1)})^{n_{\ell-1}}$.

This realization gives rise to a Cauchy formula (see~\cite{NCSF2} for the
$l=1$ case), which in turn allows one to identify the dual of $\NCSF^{(\ell)}$
with an algebra introduced by S. Poirier in~\cite{Poi}.
It is detailed in the following section.

Note that $\NCSF^{(\ell)}$ is the natural noncommutative analog of McMahon's
algebra of multisymmetric functions~\cite{McM, Ge87}.

\section{Quasi-symmetric functions of level $\ell$}


\subsection{Cauchy formula of level $\ell$}

Let $\X = X^0 \sqcup \cdots \sqcup X^{\ell-1}$,
where $X^i=\{ x_j^{(i)},j\geq1\}$, be an $\ell$-colored alphabet of
commutative variables, also commuting with $\A$.
Imitating the level $1$ case (see~\cite{NCSF6}), we define the Cauchy kernel

\begin{equation}
\label{CauchyXA}
K(\X,\A) = \prod_{j\geq1}^{\rightarrow}
\sigma_{\left(x_j^{(0)}, \ldots, x_j^{(\ell-1)}\right)} (\A).
\end{equation}

Expanding on the basis $S^{\bf I}$ of $\NCSF^{(\ell)}$, we get as coefficients
what can be called the \emph{level $\ell$ monomial quasi-symmetric functions}
$M_{\bf I}(\X)$

\begin{equation}
\label{K-expand}
K(\X,\A) = \sum_{\bf I} M_{\bf I}(\X) S^{\bf I}(\A),
\end{equation}
defined by
\begin{equation}
M_{\bf I}(\X) = \sum_{j_1<\cdots<j_m}
{\bf x}^{{\bf i}_1}_{j_1} \cdots
{\bf x}^{{\bf i}_m}_{j_m},
\end{equation}
with ${\bf I}=({\bf i}_1,\ldots,{\bf i}_m)$.

These functions form a basis of a subalgebra $\QSym^{(\ell)}$ of
$\K[\X]$, which we shall call the \emph{algebra of quasi-symmetric functions
of level $\ell$}.

\subsection{Poirier's Quasi-symmetric functions}

The functions $M_{\bf I}(\X)$ can be recognized as a basis of one of the
algebras introduced by Poirier: the $M_{\bf I}$ coincide with the $M_{(C,v)}$
defined in~\cite{Poi}, p.~324, formula (1), up to indexation.

Following Poirier, we introduce the level $\ell$ quasi-ribbon functions by
summing over an order on $\ell$-compositions:
an $\ell$-composition $C$ is finer than $C'$, and we write $C\finer C'$, if
$C'$ can be obtained by repeatedly summing up two consecutive columns of $C$
such that no nonzero element of the left one is strictly below a nonzero
element of the right one.

This order can be described in a much easier and natural way if one recodes
an $\ell$-composition ${\bf I}$ as a pair of words, the first one $d({\bf I})$
being the set of sums of the elements of the first $k$ columns of $\bf I$, the
second one $c({\bf I})$ being obtained by concatenating the words
$i^{{\bf I}_{i,j}}$ while reading of $\bf I$ by columns, from top to bottom
and from left to right.
For example, the $3$-composition of Equation~(\ref{exM}) satisfies
\begin{equation}
d({\bf I}) = \{5, 10, 14, 19\} \text{\quad and\quad}
c({\bf I}) = 13333\, 22233\, 1123\, 12333\,.
\end{equation}
Moreover, this recoding is a bijection if the two words $d({\bf I})$ and
$c({\bf I})$ are such that the descent set of $c({\bf I})$ is a subset of
$d({\bf I})$.
The order previously defined on $\ell$-compositions is in this context the
inclusion order on sets $d$: $(d',c)\finer (d,c)$ iff $d'\subseteq d$.

It allows us to define the \emph{level $\ell$ quasi-ribbon functions}
$F_{\bf I}$
by
\begin{equation}
F_{\bf I} = \sum_{{\bf I'}\finer {\bf I}} M_{\bf I'}.
\end{equation}
Notice that this last description of the order $\finer$ is reminiscent of the
order $\finer'$ on descent sets used in the context of quasi-symmetric
functions and non-commutative symmetric functions: more precisely, since it
does not modify the word $c({\bf I})$, the order $\finer$ restricted to
$\ell$-compositions of weight $n$ amounts for $\ell^n$ copies of the order
$\finer'$.
The computation of its M\"obius function is therefore straightforward.

Moreover, one can obtain the $F_{\bf I}$ as the commutative image
of certain $\F_{\sigma,u}$: any pair $(\sigma,u)$ such that $\sigma$ has
descent set $d({\bf I})$ and $u=c({\bf I})$ will do.

\subsection{Coproducts and alphabets}

\subsubsection{}
Recall that to define the $\G_{\sigma,u}(\A)$ of an $\ell$-colored alphabet
$\A=A\times C$, we only need a total order on $A$. Hence, if $B$ is another
copy of $A$ commuting with $A$, we can define $\A+\B$ as $(A+B)\times C$ where
$A+B$ is the ordinal sum, and thus make sense of $\G_{\sigma,u}(\A+\B)$.

As usual, we identify $F(\A)G(\B)$ with $F\otimes G$.

\begin{lemma}
For any $F\in\FQSym^{(\ell)}$,
$\F(\A+\B)\in \FQSym^{(\ell)}\otimes \FQSym^{(\ell)}$, and
\begin{equation}
\F(\A+\B)=\Delta F\,,
\end{equation}
where $\Delta$ is the coproduct defined by~(\ref{deltaG}).
\end{lemma}

\Proof
It is sufficient to show this for $\ell=1$, which is done in~\cite{NCSF6}.
\qed 

Let us observe that on this picture, it is clear that the restriction
of $\Delta$ to $\Sym^{(\ell)}$ is dual to the product of $QSym^{(\ell)}$.
By definition of the Cauchy kernel (\ref{CauchyXA}), we have
\begin{equation}
K(\X,\A+\B)= K(\X,\A)K(\X,\B)\,,
\end{equation}
and by~(\ref{K-expand}), this implies that $\Delta$ is dual of the
multiplication of $QSym^{(\ell)}$. 

\subsubsection{}
The same description can be applied to the quasi-symmetric side. Let
$\X=X\times C$ and $\Y=Y\times C$, where $Y$ is a copy of $X$. Again, define
$X+Y$ as the ordinal sum of $X$ and $Y$, and $\X+\Y=(X+Y)\times C$. 

\begin{lemma}
The map $\nabla\ :\  F\mapsto F(\X+\Y)$ defines a coproduct on 
$QSym^{(\ell)}$, which is dual to the product of $\Sym^{(\ell)}$.
\end{lemma}

\Proof This follows from the identity
\begin{equation}
K(\X+\Y,\A)=K(\X,\A)K(\Y,\A)\,.
\end{equation}
\qed

\subsubsection{The internal coproduct}

From now on, we assume that the color set $C$ is an additive semigroup, such
that every element $\gamma\in C$ has a finite number of decompositions
$\gamma=\alpha+\beta$.

We define the $C$-product $\T=\X\times_C \Y$ of two $C$-colored alphabets by
\begin{equation}
\T^{(\gamma)}
=\{t^{(\gamma)}_{rs}=\sum_{\alpha+\beta=\gamma}x_r^{(\alpha)}y_s^{(\beta)}\}\,,
\end{equation}
with the pairs $(r,s)$ ordered lexicographically.

\begin{proposition}
\label{prop-intpro}
The map $\delta\ :\ F\mapsto F(\X\times_C \Y)$ defines a coassociative
coproduct on $QSym^{(\ell)}$.
\end{proposition}

We define the internal product $*$ of $\Sym^{(\ell)}$ as the dual product of
the map $\delta$.

\Proof
The coassociativity condition
\begin{equation}
(\delta \otimes Id) \circ \delta = (Id\otimes\delta) \circ \delta
\end{equation}
translates as the associativity of the $C$-product
\begin{equation}
(\X \times_C \Y) \times_C \Z =
\X \times_C (\Y \times_C \Z)
\end{equation}
which is clear since both sides are by definition
\begin{equation}
\left\{ t_{pqr}^{\mu} = \sum_{\alpha+\beta+\gamma=\mu}
                   x_p^{(\alpha)} y_q^{(\beta)} z_r^{(\gamma)}
\right\}
\end{equation}
with the lexicographic order on triples $(p,q,r)$.
\qed

\begin{example}
{\rm
With $\ell=2$ and $C=\Z/2\Z$,

\scriptsize
\begin{equation*}
\begin{split}
S^{\hbox{\scriptsize$\begin{pmatrix}1&0\\1&1\end{pmatrix}$}} *
S^{\hbox{\scriptsize$\begin{pmatrix}0&2\\1&0\end{pmatrix}$}} &=
\mu_2 \left[ \left(
    S^{\hbox{\scriptsize$\begin{pmatrix}1\\1\end{pmatrix}$}}\otimes
      S^{\hbox{\scriptsize$\begin{pmatrix}0\\1\end{pmatrix}$}} \right)
    *_2
    \Delta S^{\hbox{\scriptsize$\begin{pmatrix}0&2\\1&0\end{pmatrix}$}}
    \right]\\
& =  \left(S^{\hbox{\scriptsize$\begin{pmatrix}1\\1\end{pmatrix}$}}*
           S^{\hbox{\scriptsize$\begin{pmatrix}2\\0\end{pmatrix}$}}\right)
     \left(S^{\hbox{\scriptsize$\begin{pmatrix}0\\1\end{pmatrix}$}}*
           S^{\hbox{\scriptsize$\begin{pmatrix}0\\1\end{pmatrix}$}}\right)
  +  \left(S^{\hbox{\scriptsize$\begin{pmatrix}1\\1\end{pmatrix}$}}*
           S^{\hbox{\scriptsize$\begin{pmatrix}0&1\\1&0\end{pmatrix}$}}\right)
     \left(S^{\hbox{\scriptsize$\begin{pmatrix}0\\1\end{pmatrix}$}}*
           S^{\hbox{\scriptsize$\begin{pmatrix}1\\0\end{pmatrix}$}}\right) \\
& = S^{\hbox{\scriptsize$\begin{pmatrix}1&1\\1&0\end{pmatrix}$}} 
  + S^{\hbox{\scriptsize$\begin{pmatrix}1&1&0\\0&0&1\end{pmatrix}$}} 
  + S^{\hbox{\scriptsize$\begin{pmatrix}0&0&0\\1&1&1\end{pmatrix}$}}.
\end{split}
\end{equation*}
}
\end{example}

\subsubsection{The splitting formula}
The definition of $*$ by duality with the $C$-product implies that it
satisfies the splitting formula. We just need to check a trivial property:

\begin{lemma}
The $C$-product is distributive over the colored ordinal sum:
\begin{equation}
(\X'+\X'')\times_C \Y= \X'\times_C+\X''\times_C\Y\,.
\end{equation}
\end{lemma}
\qed

\begin{proposition}
\label{split-prop}
Let $\mu_r: (\NCSF^{(\ell)})^{\otimes r} \to \NCSF^{(\ell)}$ be the product
map.
Let $\Delta^{(r)} : (\NCSF^{(\ell)}) \to (\NCSF^{(\ell)})^{\otimes r}$ be the
$r$-fold coproduct, and $*_r$ be the extension of the internal product to
$(\NCSF^{(\ell)})^{\otimes r}$.
Then, for $F_1,\ldots,F_r$, and $G\in\NCSF^{(\ell)}$,

\begin{equation}
\label{split-form}
(F_1\cdots F_r) * G = \mu_r [ (F_1\otimes\cdots\otimes F_r) *_r
\Delta^{(r)}G ].
\end{equation}
\end{proposition}

\Proof It is enough to consider the case $r=2$. Let $F\in QSym^{(\ell)}$
and $U,W,W\in\Sym^{(\ell)}$. We have, writing $\X\Y=\X\times_C \Y$ for short,
and assuming duality between $\X$ and $\A$, $\Y$ and $\B$ and so on,
\begin{equation*}
\begin{split}
\<F, (UV)*W\> &= \<F(\X'\X''), (UV)(\A')W(\A'')\>\\
&= \<F(\X'+\Y')\X''), U(\A')V(\B')W(\A'')\>\\
&= \<F(\X'\X''+\Y'\Y''),U(\A')V(\B')W(\A''+\B'')\>\\
&= \<F(\X'+\Y'), U(\A')V(\B')W(\A'+\B')\>\\
&= \<\nabla F, (U\otimes V)*_2\Delta W\>=\<F,\mu[(U\otimes V)*_2\Delta W]\>\,.
\end{split}
\end{equation*}
\qed

\subsubsection{Evaluation of internal products} 

Let us start with the simplest case, $S_{\bf i}*S_{\bf j}$.
The coefficient of $S^{\bf K}$ in this product is the coefficient of
$M_{{\bf i}}\otimes M_{\bf j}$ in $\delta M_{\bf K}$, which is also the
coefficient of ${\bf x}^{\bf i} {\bf y}^{\bf j}$ in
$M_{\bf K}({\bf x} \times_C {\bf y})$.
This is zero if $\bf K$ has more than one column, so that
\begin{equation}
\label{ijn}
S_{\bf i} * S_{\bf j} = \sum_{\bf n} d_{\bf i j}^{\bf n} S_{\bf n}
\end{equation}
contains only one-part vector compositions.

\begin{lemma}
\label{lem-dij}
If the color semigroup is $C=\N$, then the coefficient
$d_{\bf i j}^{\bf n}$ of $S_{\bf n}$ in
$S_{\bf i} * S_{\bf j}$ is equal to the coefficient of the monomial
symmetric function $m_\mu$ in the product $m_\alpha m_\beta$,
where $\mu$ is the partition
$(0^{n_0}1^{n_1}\dots)$,
$\alpha = (0^{i_0}1^{i_1}\dots)$,
$\beta = (0^{j_0}1^{j_1}\dots)$,
the monomial functions being taken over an alphabet of $n$ letters,
where $n=|{\bf i}|$.
\end{lemma}

\Proof
From~(\ref{ijn}), we only need to compute the coproducts
$\delta M_{\bf n}$. For this, it is sufficient to use alphabets of the form
${\bf x}=\{x\}\otimes C$,
${\bf y}=\{y\}\otimes C$.
Then,
\begin{equation}
M_{\bf n} ({\bf x}\times_C {\bf y}) =
\prod_{k\geq0} \ \left(\sum_{i+j=k} x^{(i)} y^{(j)}\right)^{n_k},
\end{equation}
and we see that the coefficient of
$M_{\bf i}({\bf x}) M_{\bf j}({\bf y})=
{\bf x^i y^j}$
is the same as the coefficient of $h_\alpha \otimes h_\beta$
in
\begin{equation}
\prod_{k\geq0} \left(\sum_{i+j=k} h_i\otimes h_j\right)^{n_k}
= \Delta \left(h_0^{n_0} h_1^{n_1} \dots \right) = \Delta h_\mu.
\end{equation}
\qed

For example,
\begin{equation}
S^{\hbox{\scriptsize$\begin{pmatrix}2\\2\end{pmatrix}$}} *
S^{\hbox{\scriptsize$\begin{pmatrix}3\\1\end{pmatrix}$}} =
3\,S^{\hbox{\scriptsize$\begin{pmatrix}1\\3\end{pmatrix}$}} +
S^{\hbox{\scriptsize$\begin{pmatrix}2\\1\\1\end{pmatrix}$}}.
\end{equation}
This result is compatible with the fact that
\begin{equation}
m_{11} m_{1} = 3\, m_{111} + m_{21}.
\end{equation}
As another example,
\begin{equation}
\label{ex021}
S^{\hbox{\scriptsize$\begin{pmatrix}0\\2\\1\end{pmatrix}$}} *
S^{\hbox{\scriptsize$\begin{pmatrix}1\\1\\1\end{pmatrix}$}} =
S^{\hbox{\scriptsize$\begin{pmatrix}0\\1\\1\\0\\1\end{pmatrix}$}} +
2\,S^{\hbox{\scriptsize$\begin{pmatrix}0\\1\\0\\2\end{pmatrix}$}} +
2\,S^{\hbox{\scriptsize$\begin{pmatrix}0\\0\\2\\1\end{pmatrix}$}}.
\end{equation}
One can then check that the previous result amounts to selecting the
partitions of size at most $3$ in
\begin{equation}
m_{211} m_{21} = m_{421} + 2\,m_{331} + 2\,m_{322} +\dots
\end{equation}

Together with the splitting formula~(\ref{split-form}),
this result determines all the products
$S^{{\bf I}} * S^{{\bf J}}$, since one also has
\begin{equation}
S^{{\bf i}} * S^{{\bf J}} = S^{{\bf J}} * S^{{\bf i}}
\end{equation}
thanks to the isomorphism of ordered colored alphabets
\begin{equation}
{ \bf x} \times_C {\bf Y} \simeq {\bf Y} \times_C {\bf x},
\end{equation}
where ${\bf x} = \{x\}\otimes C$.

When the color group is $\Z/\ell\Z$, the result is obtained by reduction
modulo $\ell$, \emph{e.g.}, with $l=3$,
we get from example (\ref{ex021})
\begin{equation}
\label{ex021b}
S^{\hbox{\scriptsize$\begin{pmatrix}0\\2\\1\end{pmatrix}$}} *
S^{\hbox{\scriptsize$\begin{pmatrix}1\\1\\1\end{pmatrix}$}} =
2\,S^{\hbox{\scriptsize$\begin{pmatrix}0\\2\\1\end{pmatrix}$}} +
2\,S^{\hbox{\scriptsize$\begin{pmatrix}2\\1\\0\end{pmatrix}$}} +
2\,S^{\hbox{\scriptsize$\begin{pmatrix}1\\0\\2\end{pmatrix}$}}.
\end{equation}
Note that the coefficient  of an $S^{\bf I}$ can change when computing modulo
$\ell$: for all pairs of parts $k$ and $k'$ added together, a factor
$\binom{k+k'}{k}$ appears.

\subsection{Generalized descent algebras}

In a preliminary draft of this work~\cite{NT-level},
we introduced the internal product in a different way.
Assuming that $C$ has a semigroup structure, we regard colored permutations as
elements of the wreath product $H = C\wr \SG_n$, and for
$h'$, $h''\in H$, we set
\begin{equation}
\G_{h'} *' \G_{h''} = \G_{h'' h'}
\end{equation}
(opposite law, as in the classical case of $\NCSF$).

\begin{theorem}
Let $C$ be a commutative semigroup.
\begin{enumerate}
\item $\NCSF_n^{(C)}$ is a subalgebra of $\FQSym^{(C)}$, for the operation
$*'$ defined previously.
\item The restriction of $*'$ to $\NCSF_n^{(C)}$ satisfies the splitting
formula~(\ref{split-form}).
\item $S^{{\bf i}} *' S^{{\bf j}}$ is given by Lemma~\ref{lem-dij}.
\item $S^{{\bf i}} *' S^{{\bf J}} = S^{{\bf J}} *' S^{{\bf i}} $.
\end{enumerate}
\end{theorem}

\Proof
The proofs of $(1)$ and $(2)$ can be found in~\cite{BH}.
$(3)$ and $(4)$ follow from the definition in~(\ref{defS}) and the internal
product on $\FQSym^{(C)}$ given by~(\ref{intpF}).
\qed

This provides an analogue of Solomon's descent algebra for the wreath
product $C\wr\SG_n$.
Note that the definition remains valid for $C=\Z$, so that we get a descent
algebra for the (extended) affine Weyl groups of type $A$,
$\widehat{\SG}_n=\Z\,\wr\, \SG_n$.
%

\subsection{Ordinary multi-symmetric functions}

A consequence of the above results is that the algebra $Sym^{(\ell)}$
of ordinary (commutative) multi-symmetric functions admits an internal
product. If we denote by $\underline{F}\in Sym^{(\ell)}$ the commutative
image of $F\in\Sym^{(\ell)}$, we have
\begin{equation}
\underline{F*G} =\underline{F'*G'}\quad\text{as soon as
$\underline{F}=\underline{F'}$ and $\underline{G}=\underline{G'}$.}
\end{equation}
For example,

\begin{equation}
\begin{split}
S^{\hbox{\tiny$\begin{pmatrix}0&3\\1&1\end{pmatrix}$}} *
S^{\hbox{\tiny$\begin{pmatrix}1&2\\1&1\end{pmatrix}$}}
=&\ \ \ \
 S^{\hbox{\tiny$\begin{pmatrix}0&0&2\\1&0&1\\0&1&0\end{pmatrix}$}} +
 S^{\hbox{\tiny$\begin{pmatrix}0&0&2\\0&1&1\\1&0&0\end{pmatrix}$}} +
 S^{\hbox{\tiny$\begin{pmatrix}0&1&2\\0&0&0\\1&0&1\end{pmatrix}$}} +
 S^{\hbox{\tiny$\begin{pmatrix}0&0&2\\1&1&0\\0&0&1\end{pmatrix}$}} \\
 &+
2S^{\hbox{\tiny$\begin{pmatrix}0&1&1\\0&0&2\\1&0&0\end{pmatrix}$}} +
2S^{\hbox{\tiny$\begin{pmatrix}0&0&1\\1&1&2\\0&0&0\end{pmatrix}$}}\\ &+
 S^{\hbox{\tiny$\begin{pmatrix}1&0&1\\0&1&1\\1&0&0\end{pmatrix}$}} +
 S^{\hbox{\tiny$\begin{pmatrix}1&0&1\\1&1&0\\0&0&1\end{pmatrix}$}} +
 S^{\hbox{\tiny$\begin{pmatrix}1&0&1\\1&0&1\\0&1&0\end{pmatrix}$}} +
 S^{\hbox{\tiny$\begin{pmatrix}1&0&2\\0&0&0\\1&1&0\end{pmatrix}$}}\\ &+
2S^{\hbox{\tiny$\begin{pmatrix}1&0&0\\1&1&2\\0&0&0\end{pmatrix}$}} +
2S^{\hbox{\tiny$\begin{pmatrix}0&0&1\\2&1&1\\0&0&0\end{pmatrix}$}} +
2S^{\hbox{\tiny$\begin{pmatrix}0&0&2\\2&0&0\\0&1&0\end{pmatrix}$}}.
\end{split}
\end{equation}

\begin{equation}
\begin{split}
S^{\hbox{\tiny$\begin{pmatrix}0&3\\1&1\end{pmatrix}$}} *
S^{\hbox{\tiny$\begin{pmatrix}2&1\\1&1\end{pmatrix}$}}
=&\ \ \ \
 S^{\hbox{\tiny$\begin{pmatrix}2&0&0\\1&0&1\\0&1&0\end{pmatrix}$}} +
 S^{\hbox{\tiny$\begin{pmatrix}2&0&0\\1&1&0\\0&0&1\end{pmatrix}$}} +
 S^{\hbox{\tiny$\begin{pmatrix}2&0&1\\0&0&0\\1&1&0\end{pmatrix}$}} +
 S^{\hbox{\tiny$\begin{pmatrix}2&0&0\\0&1&1\\1&0&0\end{pmatrix}$}} \\ &+
2S^{\hbox{\tiny$\begin{pmatrix}1&0&1\\2&0&0\\0&1&0\end{pmatrix}$}} +
2S^{\hbox{\tiny$\begin{pmatrix}1&0&0\\2&1&1\\0&0&0\end{pmatrix}$}} \\ &+
 S^{\hbox{\tiny$\begin{pmatrix}0&1&1\\1&0&1\\0&1&0\end{pmatrix}$}} +
 S^{\hbox{\tiny$\begin{pmatrix}0&1&1\\1&1&0\\0&0&1\end{pmatrix}$}} +
 S^{\hbox{\tiny$\begin{pmatrix}0&1&1\\0&1&1\\1&0&0\end{pmatrix}$}} +
 S^{\hbox{\tiny$\begin{pmatrix}0&2&1\\0&0&0\\1&0&1\end{pmatrix}$}} \\ &+
2S^{\hbox{\tiny$\begin{pmatrix}0&0&1\\1&2&1\\0&0&0\end{pmatrix}$}} +
2S^{\hbox{\tiny$\begin{pmatrix}0&1&0\\1&1&2\\0&0&0\end{pmatrix}$}} +
2S^{\hbox{\tiny$\begin{pmatrix}0&2&0\\0&0&2\\1&0&0\end{pmatrix}$}}.
\end{split}
\end{equation}

\begin{equation}
\begin{split}
S^{\hbox{\tiny$\begin{pmatrix}3&0\\1&1\end{pmatrix}$}} *
S^{\hbox{\tiny$\begin{pmatrix}1&2\\1&1\end{pmatrix}$}}
=&\ \ \ \
 S^{\hbox{\tiny$\begin{pmatrix}0&0&2\\1&0&1\\0&1&0\end{pmatrix}$}} +
 S^{\hbox{\tiny$\begin{pmatrix}0&0&2\\0&1&1\\1&0&0\end{pmatrix}$}} +
 S^{\hbox{\tiny$\begin{pmatrix}1&0&2\\0&0&0\\0&1&1\end{pmatrix}$}} +
 S^{\hbox{\tiny$\begin{pmatrix}0&0&2\\1&1&0\\0&0&1\end{pmatrix}$}} \\ &+
2S^{\hbox{\tiny$\begin{pmatrix}1&0&1\\0&0&2\\0&1&0\end{pmatrix}$}} +
2S^{\hbox{\tiny$\begin{pmatrix}0&0&1\\1&1&2\\0&0&0\end{pmatrix}$}} \\ &+
 S^{\hbox{\tiny$\begin{pmatrix}1&1&0\\0&1&1\\1&0&0\end{pmatrix}$}} +
 S^{\hbox{\tiny$\begin{pmatrix}1&1&0\\1&0&1\\0&1&0\end{pmatrix}$}} +
 S^{\hbox{\tiny$\begin{pmatrix}1&1&0\\1&1&0\\0&0&1\end{pmatrix}$}} +
 S^{\hbox{\tiny$\begin{pmatrix}1&2&0\\0&0&0\\1&0&1\end{pmatrix}$}} \\ &+
2S^{\hbox{\tiny$\begin{pmatrix}1&0&0\\1&2&1\\0&0&0\end{pmatrix}$}} +
2S^{\hbox{\tiny$\begin{pmatrix}0&1&0\\2&1&1\\0&0&0\end{pmatrix}$}} +
2S^{\hbox{\tiny$\begin{pmatrix}0&2&0\\2&0&0\\0&0&1\end{pmatrix}$}}.
\end{split}
\end{equation}

\begin{equation}
\begin{split}
S^{\hbox{\tiny$\begin{pmatrix}3&0\\1&1\end{pmatrix}$}} *
S^{\hbox{\tiny$\begin{pmatrix}2&1\\1&1\end{pmatrix}$}}
=&\ \ \ \
 S^{\hbox{\tiny$\begin{pmatrix}2&0&0\\1&0&1\\0&1&0\end{pmatrix}$}} +
 S^{\hbox{\tiny$\begin{pmatrix}2&0&0\\1&1&0\\0&0&1\end{pmatrix}$}} +
 S^{\hbox{\tiny$\begin{pmatrix}2&1&0\\0&0&0\\1&0&1\end{pmatrix}$}} +
 S^{\hbox{\tiny$\begin{pmatrix}2&0&0\\0&1&1\\1&0&0\end{pmatrix}$}} \\ &+
2S^{\hbox{\tiny$\begin{pmatrix}1&1&0\\2&0&0\\0&0&1\end{pmatrix}$}} +
2S^{\hbox{\tiny$\begin{pmatrix}1&0&0\\2&1&1\\0&0&0\end{pmatrix}$}} \\ &+
 S^{\hbox{\tiny$\begin{pmatrix}1&0&1\\0&1&1\\1&0&0\end{pmatrix}$}} +
 S^{\hbox{\tiny$\begin{pmatrix}1&0&1\\1&1&0\\0&0&1\end{pmatrix}$}} +
 S^{\hbox{\tiny$\begin{pmatrix}1&0&1\\1&0&1\\0&1&0\end{pmatrix}$}} +
 S^{\hbox{\tiny$\begin{pmatrix}2&0&1\\0&0&0\\0&1&1\end{pmatrix}$}} \\ &+
2S^{\hbox{\tiny$\begin{pmatrix}0&0&1\\2&1&1\\0&0&0\end{pmatrix}$}} +
2S^{\hbox{\tiny$\begin{pmatrix}1&0&0\\1&1&2\\0&0&0\end{pmatrix}$}} +
2S^{\hbox{\tiny$\begin{pmatrix}2&0&0\\0&0&2\\0&1&0\end{pmatrix}$}}.
\end{split}
\end{equation}

If one denotes by $h$ the commutative image of $S$, one easily checks that

\begin{equation}
\begin{split}
h^{\hbox{\tiny$\begin{pmatrix}3\\1\end{pmatrix}$}}
h^{\hbox{\tiny$\begin{pmatrix}0\\1\end{pmatrix}$}} *
h^{\hbox{\tiny$\begin{pmatrix}2\\1\end{pmatrix}$}}
h^{\hbox{\tiny$\begin{pmatrix}1\\1\end{pmatrix}$}}
=&\ \ \ \
 2h^{\hbox{\tiny$\begin{pmatrix}2\\1\\0\end{pmatrix}$}} 
  h^{\hbox{\tiny$\begin{pmatrix}0\\0\\1\end{pmatrix}$}} 
  h^{\hbox{\tiny$\begin{pmatrix}0\\1\\0\end{pmatrix}$}} +
 h^{\hbox{\tiny$\begin{pmatrix}2\\0\\1\end{pmatrix}$}}  
 h^{\hbox{\tiny$\begin{pmatrix}1\\0\\0\end{pmatrix}$}}  
 h^{\hbox{\tiny$\begin{pmatrix}0\\0\\1\end{pmatrix}$}} +
 h^{\hbox{\tiny$\begin{pmatrix}2\\0\\1\end{pmatrix}$}} 
 h^{\hbox{\tiny$\begin{pmatrix}0\\1\\0\end{pmatrix}$}} 
 h^{\hbox{\tiny$\begin{pmatrix}0\\1\\0\end{pmatrix}$}} \\ &+
2h^{\hbox{\tiny$\begin{pmatrix}1\\2\\0\end{pmatrix}$}}  
 h^{\hbox{\tiny$\begin{pmatrix}1\\0\\0\end{pmatrix}$}} 
 h^{\hbox{\tiny$\begin{pmatrix}0\\0\\1\end{pmatrix}$}} +
2h^{\hbox{\tiny$\begin{pmatrix}1\\2\\0\end{pmatrix}$}}
 h^{\hbox{\tiny$\begin{pmatrix}0\\1\\0\end{pmatrix}$}} 
 h^{\hbox{\tiny$\begin{pmatrix}0\\1\\0\end{pmatrix}$}} \\ &+
2h^{\hbox{\tiny$\begin{pmatrix}1\\1\\0\end{pmatrix}$}}  
 h^{\hbox{\tiny$\begin{pmatrix}1\\0\\1\end{pmatrix}$}}  
 h^{\hbox{\tiny$\begin{pmatrix}0\\1\\0\end{pmatrix}$}} +
 h^{\hbox{\tiny$\begin{pmatrix}1\\1\\0\end{pmatrix}$}}  
 h^{\hbox{\tiny$\begin{pmatrix}1\\1\\0\end{pmatrix}$}}  
 h^{\hbox{\tiny$\begin{pmatrix}0\\0\\1\end{pmatrix}$}} +
 h^{\hbox{\tiny$\begin{pmatrix}2\\0\\0\end{pmatrix}$}} 
 h^{\hbox{\tiny$\begin{pmatrix}1\\0\\1\end{pmatrix}$}} 
 h^{\hbox{\tiny$\begin{pmatrix}0\\0\\1\end{pmatrix}$}} \\ &+
4h^{\hbox{\tiny$\begin{pmatrix}1\\1\\0\end{pmatrix}$}}  
 h^{\hbox{\tiny$\begin{pmatrix}0\\2\\0\end{pmatrix}$}}  
 h^{\hbox{\tiny$\begin{pmatrix}0\\1\\0\end{pmatrix}$}} +
2h^{\hbox{\tiny$\begin{pmatrix}2\\0\\0\end{pmatrix}$}} 
 h^{\hbox{\tiny$\begin{pmatrix}0\\2\\0\end{pmatrix}$}} 
 h^{\hbox{\tiny$\begin{pmatrix}0\\0\\1\end{pmatrix}$}}.
\end{split}
\end{equation}

\section{The Mantaci-Reutenauer algebra}

\subsection{Monochromatic complete functions}

Let ${\bf e}_i$ be the canonical basis of $\N^\ell$. For $n\geq1$, let
\begin{equation}
S_n^{(i)} = S_{n\cdot {\bf e}_i} \in\NCSF^{(\ell)},
\end{equation}
be the \emph{monochromatic complete symmetric functions}.

\begin{proposition}
The $S_n^{(i)}$ generate a Hopf subalgebra $\MR^{(\ell)}$ of $\NCSF^{(\ell)}$,
which is isomorphic to the Mantaci-Reutenauer descent algebra of the wreath
products $\SG_n^{(\ell)} = (\Z/\ell\Z) \wr\SG_n$
if $C=\Z/\ell\Z$.
\end{proposition}

\Proof
$\MR^{(\ell)}$ is obviously stable by the product and coproduct coming from
$\NCSF^{(\ell)}$, hence is a Hopf subalgebra of $\NCSF^{(\ell)}$.

As a Hopf algebra, it is clearly isomorphic to the Mantaci-Reutenauer algebra,
having the same number of independent generators in each degree, and behaving
in the same way (divided powers) under the coproduct.
The isomorphism for the internal product comes from the splitting formula,
which gives explicitly
$S^{(I,u)} * S^{(J,v)}$.
\qed

Since $\MR^{(\ell)}$ has $\ell$ generators in each degree, 
its dimensions are given by
\begin{equation}
\frac{1}{1-\sum_{k\geq1} \ell t^k} =
1 + \ell \sum_{n\geq1} (\ell+1)^{n-1} t^n.
\end{equation}
The bases of $\MR^{(\ell)}$ are labelled by colored compositions, that is,
pairs formed by a composition and a color vector of the same length:
\begin{equation}
(I,u) = ((i_1,\ldots,i_m),(u_1,\ldots,u_m)).
\end{equation}


\subsection{Primitive elements of $\MR^{(\ell)}$}
\label{prim-MR}

$\MR^{(\ell)}$ being a subalgebra of a graded connected cocommutative Hopf
algebra, is itself a graded connected cocommutative Hopf algebra, so that,
thanks to the Cartier-Milnor-Moore theorem, it is the universal enveloping
algebra of $L^{(\ell)}$:

\begin{equation}
\MR^{(\ell)} = U(L^{(\ell)}),
\end{equation}
where $L^{(\ell)}$ is the Lie algebra of its primitive elements.
The following property is obtained in the same way
as Theorem~\ref{prim-libNC}.
\begin{theorem}
\label{prim-libMR}
As a graded Lie algebra, the primitive Lie algebra $L^{(\ell)}$
of $\MR^{(\ell)}$ is free over a set indexed by colored compositions.
\end{theorem}

For example, with $\ell=2$, the generating series of the dimensions of
$L^{(\ell)}$
is
\begin{equation}
1+ 2\,t + 3\,t^2 + 8\,t^3 + 18\,t^4 + 48\,t^5 + 116\,t^6 +
312\,t^7 + \dots
\end{equation}
With $\ell=3$, one finds
\begin{equation}
1+ 3\,t + 6\,t^2 + 20\,t^3 +  60\,t^4 + 204\,t^5 +  670\,t^6 +
2340\,t^7 + \dots
\end{equation}

More generally, the dimension of $L_n^{(\ell)}$ is given by the Witt
polynomials
\begin{equation}
q_n(\ell) :=
\left\{
\begin{array}{ll}
\ell &\text{\ \ if $n=1$},\\
\frac{1}{n} \sum_{d\,|\,n} \mu(d) (l+1)^{n/d} &\text{\ \ if $n>1$},
\end{array}
\right.
\end{equation}
so that, for $n\geq2$, the dimension of $L_n^{(\ell)}$ coincide with those
of a free Lie algebra on $\ell+1$ generators of degree $1$.

\subsection{Duality}

The duality is easily worked out by means of the appropriate Cauchy kernel.
The generating function of the complete functions is
\begin{equation}
\sigma^{\MR}_{\bf x}(\A) := 1 + \sum_{j=0}^{\ell-1} \sum_{n\geq1}
S_n^{(j)}\cdot (x^{(j)})^n,
\end{equation}
and the Cauchy kernel is as usual
\begin{equation}
K^{\MR}(\X,\A) = \prod_{i\geq1}^\rightarrow \sigma^{\MR}_{{\bf x}_i}(\A)
= \sum_{(I,u)} M_{(I,u)}(\X) S^{(I,u)}(\A),
\end{equation}
where $(I,u)$ runs over colored compositions
$(I,u) = ((i_1,\ldots,i_m),(u_1,\ldots,u_m))$.
The $M_{I,u}$ are called the \emph{monochromatic monomial quasi-symmetric
functions} and satisfy
\begin{equation}
M_{(I,u)}(\X) = \sum_{j_1<\cdots<j_m}
(x_{j_1}^{(u_1)})^{i_1} \cdots (x_{j_m}^{(u_m)})^{i_m}.
\end{equation}

\begin{proposition}
The $M_{(I,u)}$ span a subalgebra of $\C[X]$ which can be identified with the
graded dual of $\MR^{(\ell)}$ through the pairing
\begin{equation}
\langle M_{(I,u)}, S^{(J,v)} \rangle = \delta_{I,J} \delta_{u,v},
\end{equation}
where $\delta$ is the Kronecker symbol.
\end{proposition}

Note that this algebra can also be obtained by imposing the relations
\begin{equation}
x_i^{(p)} x_i^{(q)} = 0, \text{\ for $p\not=q$}
\end{equation}
on the variables of $\QSym^{(l)}$.

Baumann and Hohlweg~\cite{BH} have another construction of the dual of
$\MR^{(\ell)}$ (implicitly defined in~\cite{Poi}, Lemma~11).

\section{Level $\ell$ parking quasi-symmetric functions}

\subsection{Usual parking functions}

\subsubsection{The combinatorial objects}

Recall that a \emph{parking function} on $[n]=\{1,2,\ldots,n\}$ is a word
$\park=a_1a_2\cdots a_n$ of length $n$ on $[n]$ whose nondecreasing
rearrangement $\park^\uparrow=a'_1a'_2\cdots a'_n$ satisfies $a'_i\le i$ for
all $i$.
Let $\PF_n$ be the set of such words.
It is well-known that $|\PF_n|=(n+1)^{n-1}$.

One says that $\park$ has a \emph{breakpoint} at $b$ if
$|\{a_i\le b\}|=b$. The set of all breakpoints of $\park$ is denoted by
$\bp(\park)$.
Then, $\park\in \PF_n$ is said to be \emph{prime} if $\bp(\park)=\{n\}$
(see~\cite{Stan}).

Let $\PPF_n\subset\PF_n$ be the set of prime parking functions on $[n]$.
It can easily be shown that $|\PPF_n|=(n-1)^{n-1}$ (see~\cite{Stan2}).

Finally, one says that $\park$ has a \emph{match} at $b$ if
$|\{a_i< b\}|=b-1$ and $|\{a_i\leq b\}|\geq b$. The set of all matches
of $\park$ is denoted by $\match(\park)$.


\subsubsection{Algebraic structure on parking functions}

The algebra $\PQSym$ of parking functions~\cite{NT,NT-park} is very similar to
the algebra $\FQSym$ of permutations.

Since parking functions are closed under the shifted shuffle, one defines a
product on the vector space with basis $(\F_\park)$  by

\begin{equation}
\F_{\park'} \F_{\park''} = \sum_{\park\in \park'\shuffle \park''[k]} \F_\park.
\end{equation}
The coproduct on this basis is given by the parkization algorithm~\cite{NT}:
for $w=w_1w_2\cdots w_n$ on $\{1,2,\ldots\}$, let us define
\begin{equation}
\label{dw}
d(w):=\min \{i | \#\{w_j\leq i\}<i \}\,.
\end{equation}
If $d(w)=n+1$, then $w$ is a parking function and the algorithm terminates,
returning~$w$. Otherwise, let $w'$ be the word obtained by decrementing all
the elements of $w$ greater than $d(w)$. Then $\Park(w):=\Park(w')$. Since
$w'$ is smaller than $w$ in the lexicographic order, the algorithm terminates
and always returns a parking function.

\smallskip
For example, let $w=(3,5,1,1,11,8,8,2)$. Then $d(w)=6$ and the word
$w'=(3,5,1,1,10,7,7,2)$.
Then $d(w')=6$ and $w''=(3,5,1,1,9,6,6,2)$. Finally, $d(w'')=8$ and
$w'''= (3,5,1,1,8,6,6,2)$, that is a parking function.
Thus, $\Park(w)=(3,5,1,1,8,6,6,2)$.

We then have
\begin{equation}
\Delta \F_{\park} = \sum_{\gf{u,v}{u.v=\park}}
\F_{\Park(u)} \otimes \F_{\Park(v)}.
\end{equation}

\subsubsection{Duality}

Let $\G_{\park}=\F_{\park}^* \in\PQSym^*$ be the dual basis of $(\F_\park)$.
If $\langle\,,\,\rangle$ denotes the duality bracket, the product on
$\PQSym^*$ is given by
\begin{equation}
\G_{\park'} \G_{\park''} = \sum_{\park}
    \langle\, \G_{\park'}\otimes\G_{\park''}, \Delta\F_\park \,\rangle\,
    \G_\park
= \sum_{\park \in \park'\convol\park''} \G_\park\,,
\end{equation}
where the \emph{convolution} $\park'\convol\park''$ of two parking functions
is defined as
\begin{equation}
\park'\convol\park''
= \sum_{
  \gf{u,v ; \park=u\cdot v,}{\Park(u)=\park', \Park(v)=\park''}}
\park\,.
\end{equation}

By duality, one easily gets the formula for the
coproduct of $\G_\park$ as
\begin{equation}
\Delta \G_\park := \sum_{u,v ; \park\in u\ssh v}
                   {\G_{u} \otimes \G_{v}}\,.
\end{equation}

\subsection{Colored parking functions}

Let $\ell$ be an integer, representing the number of allowed colors.
A \emph{colored parking function} of level $\ell$ and size $n$ is a pair
composed of a parking function of length $n$ and a word on $[\ell]$ of length
$n$.

Since there is no restriction on the coloring, it is obvious that there are
$\ell^n (n+1)^{n-1}$ colored parking functions of level $\ell$ and size $n$.

With two colors, one finds the sequence $a_i=2(2i+2)^{i-1}$, known as~A097629
in~\cite{Slo}:
\begin{equation}
1+2\,t + 12\,t^2 + 128\,t^3 + 2000\,t^4 + 41472\,t^5 + 1075648\,t^6 + \dots
\end{equation}

Since the convolution of two parking functions contains only parking
functions, one easily builds as in~\cite{NT} an algebra
$\PQSym^{(\ell)}$ on colored parking functions:

\begin{equation}
\G_{(\park',u')} \G_{(\park'',u'')} = \sum_{\park\in \park'\convol\park''}
  \G_{(\park,u'\cdot  u'')}.
\end{equation}

We can define a coproduct  using sums of
alphabets: again, we only need a total order on $A$ to define
the colored parkization, so that taking two isomorphic copies $A'$ and
$A''$ of $A$, we define $\A'\oplus \A''$ as $(A'\oplus A'')\times C$, where
$A'\oplus A''$ denotes the ordered sum. Assuming furthermore that $\A'$ and
$\A''$ commute, we identify $f(\A') g(\A'')$ with $f\otimes g$ and define a
coproduct by:
\begin{equation}
\Delta f(\A) = f(\A'\oplus \A'').
\end{equation}

By construction, this is an algebra morphism from
$\PQSym^{(\ell)}$ to $\PQSym^{(\ell)}\otimes \PQSym^{(\ell)}$, so that

\begin{theorem}
$\PQSym^{(\ell)}$ is a graded connected bialgebra, hence a Hopf algebra.
More precisely, the coproduct can be computed in the following way:
\begin{equation}
\label{deltaGp}
\Delta\G_{(\park,u)}
=\sum_{\gf{(\park',\park'',u',u'')}%
          {(\park,u)\in(\park',u')\ssh(\park'',u'')}}
 \G_{(\park',u')} \otimes \G_{(\park'',u'')}.
\end{equation}
\end{theorem}

\Proof
Straightforward from the previous definitions.
\qed

For example,
\begin{equation}
\Delta \G_{(41142,22115)}
=  1 \otimes \G_{(41142,22115)}
 + \G_{(112,215)}\otimes \G_{(11,21)}
 + \G_{(41142,22115)} \otimes 1.
\end{equation}

\subsection{Duality}

Let $\F_{(\park,u)}=\G_{(\park,u)}^* \in{\PQSym^{(\ell)}}^*$ be the dual basis
of $(\G_\park)$.
If $\langle\,,\,\rangle$ denotes the duality bracket, the product on
${\PQSym^{(\ell)}}^*$ is given by
\begin{equation}
\F_{(\park',u')} \F_{(\park'',u'')}
= \sum_{(\park,u)\in(\park',u')\ssh(\park'',u'')} \F_{(\park,u)}\,,
\end{equation}
where the \emph{shifted shuffle} of two colored parking functions is such that
colors follow their letters.

Using the duality bracket once more, one easily gets the formula for the
coproduct of $\F_{(\park,u)}$ as
\begin{equation}
\Delta \F_{(\park,u)}
= \sum_{\gf{(p',u'),(p'',u'')}{p'p''=\park ; u'u''=u}}
  \F_{(\Park(p'),u')} \otimes \F_{(\Park(p''),u'')}.
\end{equation}

\subsection{Algebraic structure of $\PQSym^{(\ell)}$}

Recall that a word $w$ over $\N^*$ is \emph{connected} if it cannot be
written as a shifted concatenation $w=u\sconc v$, and \emph{anti-connected} if
its mirror image $\overline{w}$ is connected.

We denote by $\connp$ the set of connected parking functions, and by
$p_n:=|\connp_n|$ the number of such parking functions of size $n$.
For later reference, we recall that the generating series of $p_n$ is
Sequence~A122708 of~\cite{Slo}:
\begin{equation*}
\begin{split}
p(t) := \sum_{n\ge 1} p_n t^n
  &= 1 - \left(\sum_{n\ge 0} (n+1)^{n-1} t^n\right)^{-1}\\
  &= t + 2\,t^2 + 11\,t^3 + 92\,t^4 + 1\,014\,t^5 +
    13\,795\,t^6 + 223\,061\,t^7 + \dots
\end{split}
\end{equation*}

Let the \emph{connected colored parking functions} be the $(\park,u)$ with
$\park$ connected and $u$ arbitrary. Their generating series is given by
$p(\ell t)$.

From~\cite{NT-park}, we immediately get
\begin{proposition}
$\PQSym^{(\ell)}$ is free over the set $\F_{\sigma,u}$ (or $\G_{\sigma,u}$),
where $(\sigma,u)$ is connected.
\end{proposition}


For example, one gets the following generating series for the
algebraic generators (connected parking functions with $\ell=2$):
\begin{equation}
2\,t + 8\,t^2 + 88\,t^3 + 1\,472\,t^4 + 32\,448\,t^5 + 882\,880\,t^6 +
28\,551\,808\,t^7 \dots
\end{equation}

\subsection{Primitive elements of $\PQSym^{(\ell)}$}

Let $\LP^{(\ell)}$ be the primitive Lie algebra of $\PQSym^{(\ell)}$.
Since $\Delta$ is not cocommutative, $\PQSym^{(\ell)}$ cannot be the universal
enveloping algebra of $\LP^{(\ell)}$.
But since it is cofree, it is, according to~\cite{LRdip}, the universal
enveloping dipterous algebra of its primitive part $\LP^{(\ell)}$.

Let $\G^{\park,u}$ be the multiplicative basis defined by
$\G^{\park,u}=\G_{\park_1,u_1}\cdots\G_{\park_r,u_r}$ where
$(\park,u)=(\park_1,u_1)\bullet\cdots\bullet(\park_r,u_r)$ is the unique
maximal factorization of $(\park,u)\in\SG_n\times C^n$ into connected colored
parking functions.

\begin{proposition}
Let $\V_{\park,u}$ be the adjoint basis of $\G^{\park,u}$.
Then, the family $(\V_{\alpha,u})_{\alpha\in\connp}$ is a basis of
$\LP^{(\ell)}$.
In particular, we have $\dim\,\LP^{(\ell)}_n =\ell^n p_n$.
Moreover, $\LP^{(\ell)}$ is free.
\end{proposition}

\Proof
The first part of the statement follows from~\cite{NCSF6}.
The second part comes from the fact that $\PQSym^{(\ell)}$ is bidendriform
(Theorem~\ref{pqsl-bidend} below).
\qed

For example, since $\LP^{(\ell)}$ is free, the generating series of the degree
of its generators is (with $\ell=2$):
\begin{equation}
  \begin{split}
    1 - \prod_{n\geq1}{(1-t^n)}^{\ell^n p_n}
    &= 1-(1-t)^2(1-t^2)^8(1-t^3)^{88} \cdots\\
    &= 2\,t + 7\,t^2 + 72\,t^3 + 1\,276\,t^4 + 28\,944\,t^5
      + 805\,288\,t^6 \\
    & + 26\,462\,232\,t^7 + \dots\\
  \end{split}
\end{equation}

\noindent
and the Hilbert series of the universal enveloping algebra of
$\PQSym^{(\ell)}$ (its domain of cocommutativity) is, again with $\ell=2$:
\begin{equation}
  \begin{split}
    \prod_{n\geq1}{(1-t^n)}^{-\ell^n c_n}
    &= 1 + 2\,t + 11\,t^2 + 108\,t^3 + 1\,713\,t^4 + 36\,470\,t^5
       + 969\,919\,t^6\\
      &+ 30\,847\,464\,t^7 + \dots \\
  \end{split}
\end{equation}

\subsection{Bidendriform and tridendriform structure}

\begin{theorem}
\label{pqsl-bidend}
The algebra $\PQSym^{(\ell)}$ has a structure of bidendriform
bialgebra, hence is free as a Hopf algebra and as a dendriform 
algebra, cofree, self-dual, and its primitive Lie algebra is free.

Moreover, the totally primitive elements of $\PQSym^{(\ell)}$ are the totally
primitive elements of $\PQSym$ with any coloring.
\end{theorem}

\Proof
It has been done in~\cite{NT-park} in the case of $\FQSym$.
But since the dendriform and codendriform structure do not involve the color
alphabet $C$, the property is true in this case as well.

Also, colors do not play any role in determining whether a given
element is (totally) primitive or not, so the last statement holds.
\qed

For example, the dendriform generators of $\PQSym$ have as degree generating
series
\begin{equation}
\sum_i {\rm dgp}_i \,t^i = t + t^2 + 7\,t^3 + 66\,t^4 + 786\,t^5 +
11\,378\,t^6 + 189\,391\,t^7 + \dots
\end{equation}
so that the dendriform generators of $\PQSym^{(2)}$ have as degree generating
series $2^i\, {\rm dgp}_i$:
\begin{equation}
2\,t + 4\,t^2 + 56\,t^3 + 1\,056\,t^4 + 25\,152\,t^5 + 721\,792\,t^6 +
24\,242\,048\,t^7 + \dots
\end{equation}

\subsubsection{Tridendriform structure}

\begin{conjecture}
As in the case of $\PQSym$, we conjecture that $\PQSym^{(\ell)}$ is a free
dendriform trialgebra.
\end{conjecture}

Note that there cannot be any relations, even tridendriform relations, among
the elements $\F_{1,c}$ where $c\in C$, so that $\PQSym^{(l)}$ contains the
free tridendriform algebra $\TD^{(\ell)}$ on $\ell$ generators.

\medskip
Recall that, if $\PQSym^{(\ell)}$ is free as a tridendriform algebra,
its generators have as generating series
\begin{equation}
{\rm TD} := \frac{{\rm PQ}-1}{2{\rm PQ}^2-{\rm PQ}},
\end{equation}
where ${\rm PQ}$ is the generating series of $\PQSym^{(\ell)}$.

For example, the tridendriform generators of $\PQSym$ have as degree
generating series
\begin{equation}
\sum_i {\rm tgp}_i \,t^i= t + 5\,t^3 + 50\,t^4 + 634\,t^5 + 9\,475\,t^6 +
163\,843\,t^7 + \dots
\end{equation}
so that the tridendriform generators of $\PQSym^{(2)}$ have as degree
generating series $2^i\, {\rm tgp}_i$:
\begin{equation}
2\,t + 40\,t^3 + 800\,t^4 + 20\,288\,t^5 + 606\,400\,t^6
+ 20\,971\,904\,t^7 + \dots
\end{equation}

\section{Type $B$ algebras}

\subsection{Parking functions of type $B$}

In~\cite{Rei}, Reiner defined non-crossing partitions of type $B$ by
analogy to the classical case. In our context, he defined the level $2$ case.
This allowed him to derive, by analogy with a simple representation
theoretical result, a definition of parking functions of type $B$ as the words
on $[n]$ of size $n$.

We shall build another set of words, also enumerated by $n^n$ that sheds light
on the relation between type-$A$ and type-$B$ parking functions and provides a
natural Hopf algebra structure on the latter.

First, fix two colors $0<1$. We say that a pair of words $(\park,u)$ composed
of a parking function and a binary colored word is a
\emph{level $2$ parking function} if
\begin{itemize}
\item the only elements of $\park$ that can have color $1$ are the matches of
$\park$.
\item for all element of $\park$ of color $1$, all letters equal to it and to
its left also have color $1$,
\item all elements of $\park$ have at least once the color $0$.
\end{itemize}

For example, there are $27$ level $2$ parking functions of size $3$: there are
the $16$ usual ones all with full color $0$, and the eleven new elements
\begin{equation}
\begin{split}
& (111,100), (111,110), (112,100), (121,100), (211,010),\\
& (113,100), (131,100), (311,010),
(122,010), (212,100), (221,100). \\
\end{split}
\end{equation}

The first time the first rule applies is with $n=4$, where one has to discard
the words $(1122,0010)$ and $(1122,1010)$ since $2$ is not a match of $1122$.
On the other hand, both words $(1133,0010)$ and $(1133,1010)$ are
$B_4$-parking functions since $1$ and $3$ are matches of $1133$.

\begin{theorem}
The restriction of $\PQSym^{(2)}$ to the $\F$ elements indexed by level $2$
parking functions is a subalgebra of $\PQSym^{(2)}$.
The restriction of $\PQSym^{(2)}$ to the $\G$ elements indexed by level $2$
parking functions is a subcoalgebra of $\PQSym^{(2)}$.
\end{theorem}

\Proof
The shifted shuffle of two $\F$ elements indexed by level $2$ parking
functions only consists in level $2$ parking functions since the definition
only involves matches (preserved by shifted shuffle) and positions of colors
$0$ and $1$ on a given letter inside a word, also preserved by shifted
shuffle.
The same property holds for the coproduct on the $\G$ side: a match on either
side of the tensor product comes from a match of the original word and all
equal letters go to the same side of the tensor product in the same order.
\qed

\subsection{Non-crossing partitions of type $B$}


Remark that in the level $1$ case, non-crossing partitions are in bijection
with non-decreasing parking functions.
To extend this correspondence to type $B$, let us start with a non-decreasing
parking function ${\bf b}$ (with no color). We factor it into the maximal
shifted concatenation of prime non-decreasing parking functions, and we choose
a color, here $0$ or $1$, for each factor.
We obtain in this way $\binom{2n}{n}$ words $\pi$, which can be identified
with \emph{type $B$ non-crossing partitions}.

Let
\begin{equation}
{\bf P}^{\pi}=\sum_{{\bf a}^\uparrow = \pi }\F_{\bf a}\,
\end{equation}
where $w^\uparrow$ denotes the nondecreasing rearrangement of the letters of
$w$. Then,
\begin{theorem}
The ${\bf P}^{\pi}$, where $\pi$ runs over the above set of non-decreasing
signed parking functions, form the basis of a cocommutative Hopf subalgebra
$\NCPQSym^{(2)}$ of $\PQSym^{(2)}$.
\end{theorem}

\Proof
The subalgebra part comes from the fact that the shifted shuffle does not mix
prime factors. The coalgebra part selects pieces of each factor, hence
satisfies that each letter of a (new) factor has identical color.
The cocommutativity part comes from the fact that all rearrangements of a
given word are considered at the same time.
\qed

All this can be extended to higher levels in a straightforward way: allow each
prime non-decreasing parking function to choose any color among $\ell$ and use
the factorization as above. Since non-decreasing parking functions are in
bijection with Dyck words, the choice can be described as: each block of a
Dyck word with no return-to-zero, chooses one color among $\ell$. In this
version, the generating series is obviously given by
\begin{equation}
\frac{1}{1- \ell\frac{1-\sqrt{1-4t}}{2}}.
\end{equation}
For $\ell=3$, we obtain Sequence~A007854 of~\cite{Slo}.

\section{Colored analogs of planar binary trees: $\PBT^{(\ell)}$}


\subsection{Definition of $\PBT^{(\ell)}$}

In the case with one color, the Hopf algebra $\PBT$ of Planar binary trees
initially defined by Loday and Ronco~\cite{LR1} can be embedded in
$\FQSym$ in the following way~\cite{HNT,HNT2}:
\begin{equation}
\label{P2F}
\Pp_T = \sum_{\sigma ; \shape(\pp(\sigma))=T} \F_\sigma\,,
\end{equation}
where $\pp$ is a simple algorithm: it is the well-known \emph{binary search
tree insertion}, such as presented, for example, by Knuth in~\cite{Kn}.

In the colored case, the definition is almost the same:
\begin{equation}
\label{P2Fcol}
\Pp_{T,u} = \sum_{(\sigma,u) ; \shape(\pp(\sigma))=T} \F_{(\sigma,u)}\,,
\end{equation}
where $u$ is a color word whose length is equal to the number of leaves of
$T$.
Note that this algebra (without its realization on words) has been previously
studied by Maria Ronco~\cite{Ron}.

Given the definitions, the generating series of the dimensions of
$\PBT^{(\ell)}$ is
\begin{equation}
1 + \ell\,t + 2\ell^2\,t^2 + 5\ell^3\,t^3 +14\ell^4\,t^4 + \dots
\end{equation}
that is, the generating series of Catalan numbers multiplied by $\ell^n$.

\subsection{Algebraic structure of $\PBT^{(\ell)}$}

Since $\PBT$ is generated by the trees with no right branch (starting from the
root), the same holds in $\PBT^{(\ell)}$:

\begin{proposition}
$\PBT^{(\ell)}$ is free over the set $\Pp_{T,u}$,
where $T$ is a tree with no right branch.
\end{proposition}

For example, the generating series of the algebraic generators of
$\PBT^{(\ell)}$ is
\begin{equation}
\ell\,t + \ell^2\,t^2 + 2\ell^3\,t^3 + 5\ell^4\,t^4 + 14\ell^5\,t^5 + \dots
\end{equation}
that is, the generating series of shifted Catalan numbers $C_{n-1}$ multiplied
by $\ell^n$.

\subsubsection{Primitive elements of $\PBT^{(\ell)}$}

Let $\L^{(\ell)}$ be the primitive Lie algebra of the algebra $\PBT^{(\ell)}$.
Since $\Delta$ is not cocommutative, $\PBT^{(\ell)}$ cannot be the universal
enveloping algebra of $\L^{(\ell)}$.
But since it is cofree, it is, according to~\cite{LRdip}, the universal
enveloping dipterous algebra of its primitive part $\L^{(\ell)}$.

Using the same arguments as in the case of $\FQSym^{(\ell)}$, one then proves

\begin{proposition}
The Lie algebra $\L^{(\ell)}$ is free.
Moreover
\begin{equation}
\dim\,\L^{(\ell)}_n =\ell^n C_{n-1},
\end{equation}
\end{proposition}

For example, since $\L^{(\ell)}$ is free, the generating series of the degree
of its generators is (with $\ell=2$):
\begin{equation}
  \begin{split}
    1 - \prod_{n\geq1}{(1-t^n)}^{\ell^n C_{n-1}}
    &= 1-(1-t)^2(1-t^2)^4(1-t^3)^{16}(1-t^4)^{80} \cdots\\
    &= 2\,t + 3\,t^2 + 8\,t^3 + 46\,t^4 + 252\,t^5\\
    & + 1\,558\,t^6 + 9\,800\,t^7 + \dots\\
  \end{split}
\end{equation}

\noindent
and the Hilbert series of the universal enveloping algebra of 
$\PBT^{(\ell)}$ (its domain of cocommutativity) is, again with $\ell=2$,
\begin{equation}
  \begin{split}
    \prod_{n\geq1}{(1-t^n)}^{-\ell^n C_{n-1}}
    &= 1 + 2\,t + 7\,t^2 + 28\,t^3 + 139\,t^4 + 762\,t^5 + 4\,549\,t^6\\
      &+ 28\,464\,t^7 + \dots \\
  \end{split}
\end{equation}

\subsubsection{Dendriform structure of $\PBT^{(\ell)}$}

Recall that $\PBT$ is the free dendriform algebra on one generator.
Since colors do not play any role in determining if a given element
is (totally) primitive or not, the same holds for $\PBT^{(\ell)}$:

\begin{proposition}
\label{pbt-bidend}
The algebra $\PBT^{(\ell)}$ is the free dendriform algebra on $\ell$
generators. It has also the structure of bidendriform bialgebra.
\end{proposition}

\begin{note}
{\rm
It is also possible to define a colored analog of $\CQSym$ the Catalan
Quasi-symmetric algebra, but the natural definition leads to a
non-cocommutative algebra, hence not sharing the basic property of $\CQSym$
itself.
}
\end{note}

\section{Examples}

\subsection{Multigraded coinvariants and colored Klyachko idempotents}

One of the very first applications of the theory of noncommutative symmetric
functions was to provide an explanation for the following coincidence.
On the one hand, consider the representation of $\SG_n$ in the coinvariant
algebra
\begin{equation}
{\mathcal H}_n = \C[x_1,\dots,x_n]_{\SG_n} = \C[x_1,\dots,x_n]/{\mathcal J},
\end{equation}
where $\mathcal J$ is the ideal generated by symmetric polynomials without
constant term. It is known~\cite{Mcd} that the graded Frobenius characteristic
of the action of $\SG_n$ on ${\mathcal H}_n$ is
\begin{equation}
\label{chqH}
\begin{split}
{\rm ch}_q({\mathcal H}_n)
&= \sum_k q^k {\rm ch}({\mathcal H}_n)^{(k)} \\
&= (q)_n h_n \left(\frac{X}{1-q}\right)
= \sum_{I\vDash n} q^{\maj(I)} r_I(X),
\end{split}
\end{equation}
where $r_I(X)$ are the ribbon Schur functions.

On the other hand, Klyachko~\cite{Kly} introduced a remarkable Lie idempotent
in $\C\SG_n$
\begin{equation}
\label{Kly}
\kappa_n = \sum_{\sigma\in \SG_n} \zeta^{\maj(\sigma)} \sigma
\end{equation}
where $\zeta$ is a primitive $n$-th root of unity.

In terms of noncommutative symmetric functions, both expressions can be
interpreted as specializations of
\begin{equation}
K_n(q) := (q)_n S_n \left(\frac{A}{1-q}\right)
 = \sum_{I\vDash n} q^{\maj(I)} R_I.
\end{equation}

This is the graded noncommutative characteristic of an action
of $H_n(0)$ on coinvariants.
This is a projective module, hence also an $\SG_n$-module,
and taking the commutative image, $(A=X)$, one
obtains the characteristic of $\SG_n$. But $K_n(q)$ can also be interpreted as
an element of the descent algebra of $\SG_n$.
A simple computation shows that for $q=\zeta$, a primitive $n$-th root of
unity, $K_n(\zeta)$ is a primitive element of $\NCSF$ of commutative image
$p_n/n$, hence is a Lie idempotent (see~\cite{NCSF2}).
Actually, $K_n(q)$ is the noncommutative Hall-Littlewood function
$\tilde{H}_{1^n}(A;q)$, and this specialization property is a special case of
a noncommutative version~\cite{Hiv} of the classical
property~\cite{Hiv,LLT,Mcd}.

There is a similar phenomenon here. Let
$q_1,\dots,q_n$ be independent variables, and consider the noncommutative
symmetric function
\begin{equation}
\KK_n(A;q_1,\dots,q_n) :=
\sum_{I\vDash n} \qq^{\MAJ(I)} R_I,
\end{equation}
where
\begin{equation}
\qq^{\MAJ(i_1,\dots,i_r)}
  := (q_1\dots q_{i_1})^r (q_{i_1+1}\dots q_{i_1+i_2})^{r-1}
     \dots (q_{i_1+\dots+i_{r-2}}\dots q_{i_1+\dots+i_{r-1}}).
\end{equation}

For example,
\begin{equation}
\KK_3 = R_3 + q_1q_2\, R_{21} + q_1\, R_{12} + q_1^2q_2\, R_{111}.
\end{equation}
\begin{equation}
\begin{split}
\KK_4 =&\ \ \ R_4 + q_1q_2q_3\, R_{31}
            + q_1q_2\, R_{22} + q_1^2q_2^2q_3\, R_{211}\\
       &+ q_1\, R_{13} + q_1^2q_2q_3\, R_{121} + q_1^2q_2\, R_{112}
        + q_1^3q_2^2q_3\, R_{1111}.
\end{split}
\end{equation}
Its commutative image is the multigraded characteristic of ${\mathcal H}_n$
with respect to the \emph{partition degree} (\emph{cf.}~\cite{Cas}).

One may also regard $\KK_n$ as an element of
\begin{equation}
\NCSF_n^{(\Z)} \subset \FQSym_n^{(\Z)}
\end{equation}
by means of the identification
\begin{equation}
R_I(A) q_1^{\alpha_1} \dots q_n^{\alpha_n}
 = \sum_{\Des(\sigma)=I} \G_{\sigma,\alpha_1\dots\alpha_n}
\end{equation}
and writing it as
\begin{equation}
\label{kkqq}
\KK_n(A;q_1,\dots,q_n) =
\sum_{I\vDash n} R_I(A) \qq^{\MAJ(I)}
\end{equation}
Note that this element lives in the (descent) algebra of the extended affine
Weyl group of type $A$
\begin{equation}
\widehat{\SG}_n = \Z^n \ltimes \SG_n =  P \ltimes \SG_n,
\end{equation}
where $P$ is the weight lattice. One can also interpret it as an element of
the usual affine Weyl group
\begin{equation}
\widetilde{\SG}_n = Q \ltimes \SG_n
\end{equation}
where $Q$ is the root lattice
\begin{equation}
Q = \{ \alpha\in P | \sum_{i=1}^n\alpha_i =0\}.
\end{equation}

\noindent
This amounts to imposing the relation
\begin{equation}
\label{prodq}
q_1\dots q_n =1,
\end{equation}
which replaces the root of unity condition $q^n=1$ in the one-parameter case.

It has been proved by McNamara and Reutenauer~\cite{mNR} that under
condition~(\ref{prodq}), $\KK_n(A;q_1,\dots,q_n)$ was a Lie idempotent in
$\C \widetilde{\SG}_n$.
Within the formalism of the present paper, this can be seen as follows: the
authors of~\cite{mNR} introduce a twisted product on
$\AA_n=\C(x_1,\ldots,x_n)[\SG_n]$ by the formula
\begin{equation}
  f({\bf x})\sigma\cdot g({\bf x})\tau
= f({\bf x})\sigma[g({\bf x})]\sigma\tau\,,
\end{equation}
where permutations act on functions as automorphisms, \emph{i.e.},
$\sigma(x_i)=x_{\sigma(i)}$, and in particular, on monomials
by $\sigma[{\bf x}^c]={\bf x}^{c\sigma^{-1}}$. Hence,
\begin{equation}\label{twpr}
(\sigma {\bf x}^c)\cdot(\tau {\bf x}^d)=\sigma\tau {\bf x}^{c\tau+d}\,,
\end{equation} 
which is the same as Formula~(\ref{wreath}) so that $\AA_n$ can be identified
with the homogeneous component of degree $n$ of $\FQSym^{(\Z)}$, by setting
\begin{equation}
\label{idmnr}
\sigma{\bf x}^c\equiv \G_\sigma {\bf x}^c \equiv \G_{\sigma,c}\,.
\end{equation}
McNamara and Reutenauer then introduce the formal series in $\FQSym^{(\Z)}$
\begin{equation}
\Theta({\bf x})=
\sum_{n\ge 0} \sum_{\sigma\in\SG_n}
\frac{\prod_{j\in\Des(\sigma)}x_{\sigma(1)}\cdots x_{\sigma(j)}}
     {\prod_{i=1^n}(1-x_{\sigma(1)}\cdots x_{\sigma(i)})}
\sigma
\end{equation}
which, applying (\ref{twpr}), and under the identification (\ref{idmnr}), 
can be rewritten as
\begin{equation}
\Theta({\bf x})=\prod_{l\ge 0}^{\leftarrow}\sum_{n\ge 0}\G_{\id_n,l^n}
= \prod_{l\ge 0}^{\leftarrow}\sigma_1(A^{(l)})\,.
\end{equation}
Indeed, introducing the new variables
\begin{equation}
y_j =x_1x_2\cdots x_j\,,
\end{equation}
and applying (\ref{twpr}) we can write
\begin{equation}
\begin{split}
\Theta({\bf x})&=
\sum_{n\ge 0} \sum_{\sigma\in\SG_n}
\G_\sigma
\frac{\prod_{d\in\Des(\sigma)}y_d}{(1-y_1)(1-y_2)\cdots (1-y_n)}\\
&=\sum_{n\ge 0}\KK_n(y_1,\ldots,y_{n-1})\frac1{((y))_n}\,,
\end{split}
\end{equation}
where $((y))_n=(1-y_1)(1-y_2)\cdots (1-y_n)$ and
\begin{equation}
\KK_n(y_1,\ldots,y_{n-1})=\sum_{\sigma\in\SG_n}
\G_\sigma \prod_{d\in\Des(\sigma)}y_d
\end{equation}
is the twisted version of the multiparameter Klyachko element introduced
in \cite[(11)]{NCSF2}. By Moebius inversion over the lattice of compositions
of $n$, we have
\begin{equation}
\KK_n(y_1,\ldots,y_{n-1}) 
=\sum_{I\vDash n}R_I\cdot
\prod_{d\in\Des(I)}y_d\\
\end{equation}
\begin{equation}
\KK_n(y_1,\ldots,y_{n-1}) \frac{1}{((y))_n} 
=\sum_{J\vDash n}
     S^J\cdot \frac{1}{1-y_n}\prod_{d\in\Des(J)}\frac{y_d}{1-y_d}\,,
\end{equation}
which implies the expression 
\begin{equation}
\Theta({\bf x})
=\prod_{l\ge 0}^{\leftarrow}\left( \sum_{n\ge 0}S_ny_n^l\right)\,.
\end{equation}
Each factor of the right-hand side is grouplike (for the coproduct
of $\FQSym^{(l)}$):
\begin{equation}
\Delta \sum_{n\ge 0}S_ny_n^l
=\sum_{n\ge 0}\sum_{i+j=n}S_i y_i^l\otimes S_jy_j^l
=\sum_{i\ge 0}S_i y_i \otimes \sum_{j\ge 0}S_jy_j^l
\end{equation}
so that also
\begin{equation}
\Delta\Theta({\bf x}) =\Theta({\bf x})\otimes\Theta({\bf x})\,.
\end{equation}
Extracting the term of degree $n$, we find
\begin{equation}
\Delta\KK_n = \sum_{i+j=n}\left(
\KK_i\frac1{((y))_i}\otimes\KK_j\frac1{((y))_j}
\right)
((y))_n
\end{equation}
so that if we send $y_n$ to 1, all terms vanish except the extreme ones,
and we get a primitive element. This is the main result of \cite{mNR}.

\subsection{A formula of Raney}

Raney~\cite{Ran} gave a combinatorial interpretation of the coefficients of
the unique solution $g(t)\in \Q[Y,Z][[t]]$ of the functional equation
\begin{equation}
g(t) = t \sum_{k=1}^\ell y_k e^{z_k g(t)},
\end{equation}
with $g(t) = \sum_{n\geq0} \frac{g_n}{n!} t^{n+1}$.
This defining equation is of the form
\begin{equation}
g(t) = t \phi(g(t)),
\end{equation}
with $\phi(u) = \sum_{k=1}^\ell y_k e^{z_k u}$.
Hence, the coefficient $g_n$ of $t^{n+1}$ in $g(t)$ is
\begin{equation}
\begin{split}
g_n &= \frac{1}{n+1}[u^n]\phi(u)^{n+1)} \\
&= \frac{1}{n+1} \sum_{\gf{n_1+\dots+n_\ell=n+1}{q_1+\dots+q_\ell=n}}
   \begin{pmatrix} n+1\\ n_1,\dots,n_\ell \end{pmatrix}
   \begin{pmatrix} n\\ q_1,\dots,q_\ell \end{pmatrix}
   y_1^{n_1}\dots y_\ell^{n_\ell} (n_1z_1)^{q_1}\dots (n_\ell z_\ell)^{q_\ell}.
\end{split}
\end{equation}
Thus, $g_n\in \N[Y,Z]$. Its combinatorial interpretation can be mechanically
derived by means of a colored version of the noncommutative Lagrange inversion
formula as formulated in~\cite{NTlag,Ges}. Consider the functional equation
\begin{equation}
\label{lag-col}
g = \sum_{k,n} b_k S_n^{(k)} g^n,
\end{equation}
where $S_n^{(k)}=S_n(A^{(k)})$ is a colored complete function and $b_k$ are
noncommuting letters. We can set
\begin{equation}
S_n = \sum_k  b_k S_n^{(k)},
\end{equation}
so that~(\ref{lag-col}) can be rewritten as
\begin{equation}
g = \sum_{n\geq0} S_n g^n,
\end{equation}
\noindent
and the solution of~\cite{NTlag} reads
\begin{equation}
\begin{split}
g &= S^0 + S^{10} + (S^{200}+S^{110}) + \dots \\
&= \sum_{\pi\in {\rm NDPF}} S^{\Ev(\pi).0},
\end{split}
\end{equation}
where {\rm NDPF} is the set of nondecreasing parking functions. Note that
$S^0=\sum b_k$ is \emph{a priori} different from $1$, and does not commute
with the other $S^i$.
Each term $S^{\Ev(\pi).0}$ represents an ordered tree $T$ in Polish notation,
so that, for example
\entrymodifiers={+<4pt>}
\begin{equation}
\vcenter{\xymatrix@C=2mm@R=2mm{
*{} & {\circ}\ar@{-}[dl]\ar@{-}[dr] \\
{\circ} & *{}                         & {\circ}\ar@{-}[d] \\
*{} & *{}                         & {\circ}\\
      }}
\end{equation}
is $S^{2010}$.

Replacing each $S^i$ by $\sum_{k=1}^\ell b_k S_i^{(k)}$, the expression
$S^{\Ev(\pi).0}$ becomes a sum over all $\ell$-colorings of the tree $T$, so
that, one recovers the combinatorial interpretation of Raney (proved in a
different way by Foata~\cite{Foa}):
let
${\bf n} =(n_1,\dots,n_\ell)$ and ${\bf q}=(q_1,\dots,q_\ell)\in \N^\ell$ and
let
$B({\bf n},{\bf q})$ be the set of $\ell$-colored trees on $n=|{\bf n}|$
vertices with $n_k$ vertices of color $k$ and such that the sum of the
arities of vertices of color $k$ is $q_k$.
Then
\begin{equation}
g = \sum_{n\geq0} \frac{1}{n!} \sum |B({\bf n},{\bf q})|
    {\bf y}^{\bf n} {\bf z}^{\bf q}.
\end{equation}


\end{document}